\documentclass{amsart}
\usepackage{amsmath,amssymb, amscd, stmaryrd}
\usepackage{fullpage}
\usepackage{enumerate}
\usepackage{cite}
\usepackage{framed}
\usepackage[dvips,dvipdf]{graphicx}
\usepackage[usenames,dvipsnames]{color}
\usepackage{color}
\usepackage{tabularx}
\usepackage{array}
\usepackage{multirow}

\newtheorem{theorem}{Theorem}[section]
\newtheorem{proposition}[theorem]{Proposition}

\theoremstyle{definition}
\newtheorem{definition}[theorem]{Definition}
\theoremstyle{definition}

\theoremstyle{remark}
\newtheorem{remark}[theorem]{Remark}
\newtheorem{example}[theorem]{Example}

\newcommand{\mb}[1]{\ensuremath{\mathbf{#1}}}

\newcommand{\NN}{\mathbb N}
\newcommand{\PP}{\mathbb P}
\newcommand{\QQ}{\mathbb Q}
\newcommand{\RR}{\mathbb R}

\newcommand{\cB}{\mathcal B}
\newcommand{\cH}{\mathcal H}
\newcommand{\cT}{\mathcal T}
\newcommand{\cV}{\mathcal V}
\newcommand{\cE}{\mathcal E}
\newcommand{\cF}{\mathcal F}
\newcommand{\cI}{\mathfrak I}

\definecolor{Blue}{rgb}{0.3,0.3,0.9}
\definecolor{Orange}{rgb}{1,0.5,0}
\definecolor{Green}{rgb}{0,0.6,0}

\author[A. Anand]{Akash Anand} \address{Akash Anand, Department of
  Mathematics and Statistics, Indian Institute of Technology, Kanpur, UP 208016}
\email{akasha@iitk.ac.in}

\author[J. Ovall]{Jeffrey S. Ovall} \address{Jeffrey S. Ovall,
  Fariborz Maseeh Department of Mathematics and Statistics, Portland
  State University, Portland, OR 97201}
\email{jovall@pdx.edu}

\author[S. Wei\ss er]{Steffen Wei\ss er} \address{Steffen Wei\ss er,
  Department of Mathematics,
  Saarland University, 66041 Saarbr\"ucken, Germany}
\email{weisser@num.uni-sb.de}

\begin{document}
\title[Nystr\"om FEM]{A Nystr\"om-based Finite Element Method on Polygonal Elements} \date{\today}

\begin{abstract}
We consider families of finite elements on polygonal meshes, that are defined implicitly on
each mesh cell as solutions of local Poisson problems with polynomial
data.  Functions in the local space on each mesh cell are evaluated via Nystr\"om
discretizations of associated integral equations, allowing for
curvilinear polygons and non-polynomial boundary data.  Several experiments
demonstrate the approximation quality of interpolated functions in
these spaces.
\end{abstract}

\maketitle

{\small

}

\section{Introduction}\label{Introduction}
During the past several years there has been increasing interest in
developing flexible finite element discretization schemes for use on
polygonal and polyhedral meshes.  Some of the appeal of such meshes is
due to the fact that refinement and coarsening, which are essential
components of high-performance computing, are much simpler when one is
not restricted to a small class of element shapes (e.g. triangles,
quadrilaterals, tetrahedra, etc.) and does not have to deal with
``hanging nodes''---allowing two edges of a polygon to meet at a
straight angle removes notion of hanging nodes altogether.  Virtual
Element Methods (VEM)
(cf.~\cite{vemMMMAS2013,vemCMA2013,vemIMA2014,vemMMMAS2014,vemIMAJNA2014,vemSINUM2014,vemCMAME2014,MR3509090,MR3564679,ANTONIETTI2017}),
which have drawn inspiration from mimetic finite difference schemes,
constitute one active line of research in this direction.  Another
involves Boundary Element-Based Finite Element Methods (BEM-FEM)
(cf.~\cite{bemfemLNCSE2009,bemfemETNA2010,bemfemJNM2011,bemfemNM2011,bemfemSINUM2012,bemfemJCAM2014,bemfemCMA2014,HofreitherLangerWeisser2016,Weisser2017}),
which have looked more toward the older Trefftz methods for
motivation.  For simple diffusion problems BEM-FEM is related to VEM
in the sense that, in most their
basic forms (cf. ~\cite{vemMMMAS2013,bemfemLNCSE2009}), both approaches arrive at the same local and global
finite element spaces in their derivations, whose functions are
described implicitly by local Poisson problems.  
An important practical difference between VEM
and BEM-FEM is how these implicit spaces are used in the formation of
stiffness (and mass) matrices, and these naturally lead to differences
in the theoretical development as well.  A third line of research
involves the development of generalized barycentric coordinates
(cf.~\cite{GRB2014,FGS2013,RGB2011a,GRB2010,NRS2014} and the references in~\cite{Floater2015}), in which
explicit bases are constructed that mimic certain key properties of
standard barycentric coordinates.  These three approaches typically
yield globally-conforming discretizations, but there has also been
significant recent activity in the development of various
non-conforming methods for polyhedral meshes.  We mention Compatible
Discrete Operator (CDO),  Hybrid High-Order (HHO) schemes
(cf.~\cite{Bonelle2014,Bonelle2015a,Bonelle2015b,DiPietro2014,DiPietro2015a,DiPietro2015b}),
Weak Galerkin (WG) schemes
(cf.~\cite{WangYe2013,WangYe2014,WangWang2014,MuWangYe2015c,MuWangYe2015b,MuWangYe2015,WangYe2016})
and discontinuous Galerkin (hp-DG) schemes (cf.~\cite{Houston2016a,Houston2016b})
in this regard.

The present work is most closely related to the BEM-FEM approach for
second-order, linear, elliptic boundary value problems posed on
polygonal domains $\Omega\subset\RR^2$: Find $u\in\cH$  
\begin{align}\label{ModelProblem}
\int_\Omega A\nabla u\cdot\nabla v+cuv\,dx=\int_\Omega f v\,dx\mbox{
for all }v\in\cH~,
\end{align}
where $\cH$ is some appropriate subspace of $H^1(\Omega)$
incorporating homogeneous Dirichlet boundary conditions, and standard
assumptions on the data $A,c,f$ ensure that the problem is coercive,
and thus well-posed.
Starting from the same
implicitly-described local spaces, we use Nystr\"om
discretizations of associated second-kind integral equations in our
evaluation of basis functions and their derivatives in the formation
of our finite element linear systems.  In contrast, BEM-FEM employs
first-kind integral equations discretized via boundary element methods
for the same purpose.  We believe that the Nystr\"om approach offers
several advantages over its boundary element counterpart in this
context, including greater ease in setting up and solving the integral
equations for higher-order discretizations, better resolution of
singular behavior in the local spaces, and the flexibility to
truly allow for elements with curved edges without modification of the core
computational kernels.  The focus of this paper is on polygonal
meshes, but we do provide some empirical insight into the behavior of
interpolation in these spaces on curved elements as well.
We finally mention the contribution~\cite{Kirsch1994}, which also employs both
finite element and Nystr\"om methods for acoustic scattering problems,
but in a very different way than that proposed here.  In that work,
finite elements are used within the scatterer, and are coupled with a
Nystr\"om approach that is employed outside the scatterer.

The paper is organized as follows: In Section~\ref{Trefftz}, we
introduce the local and global discrete approximation spaces, and
indicate how functions in the local space on a mesh cell $K$ can be
expressed implicitly in terms of solutions of integral equations posed
on $\partial K$.  The solution of such integral equations via
Nystr\"om approximations is the topic Section~\ref{Nystrom}.  An
interpolation operator is described in Section~\ref{Interpolation},
and numerical experiments demonstrate the interpolation properties on
different types of polygonal meshes. Finally, in
Section~\ref{Curvilinear}, we discuss the treatment of Dirchlet
boundary conditions, allowing for elements having curved edges along
the boundary.  In this section, we also suggest how one might allow
for elements with curved edges more generally.

\section{``Poisson Spaces'' and Associated Integral Equations}\label{Trefftz}

\begin{figure}
\begin{center} 
\includegraphics[width=1.8in,angle=90]{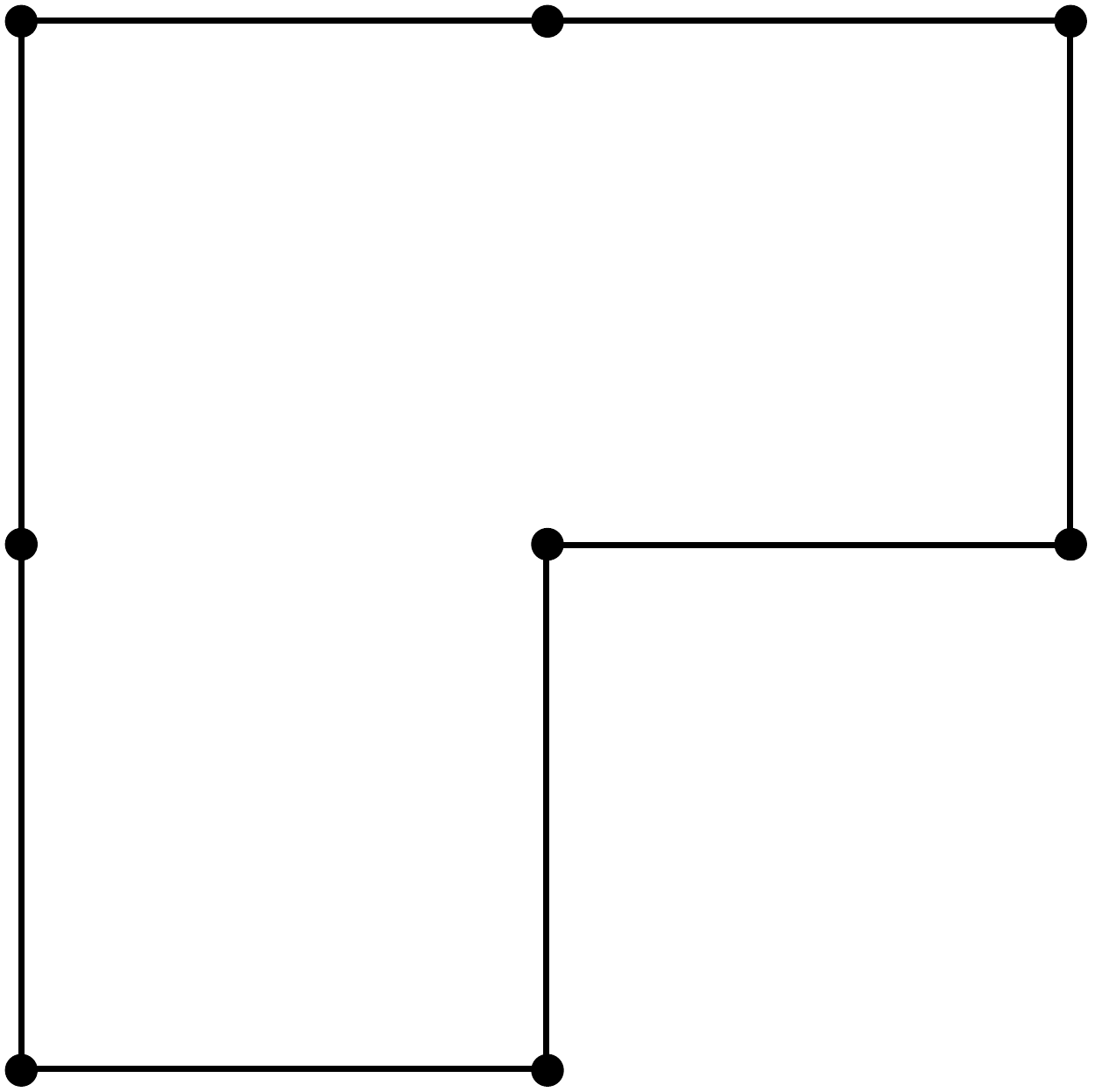} 
\includegraphics[width=1.8in]{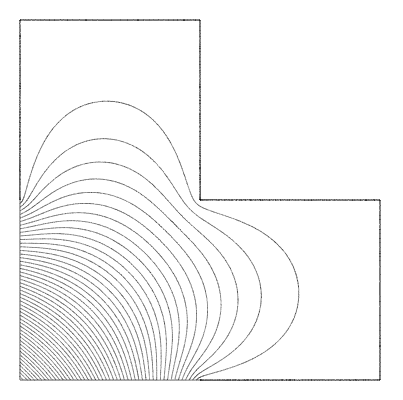} 
\includegraphics[width=1.8in]{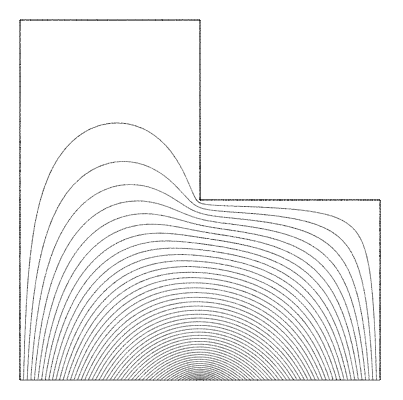} \\
\includegraphics[width=1.8in]{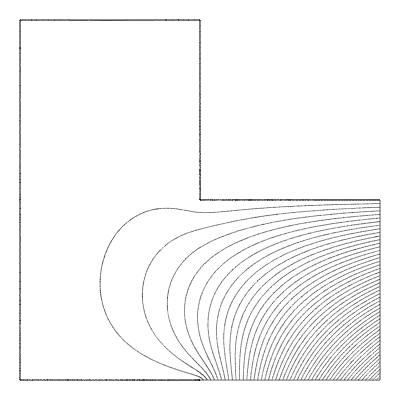} 
\includegraphics[width=1.8in]{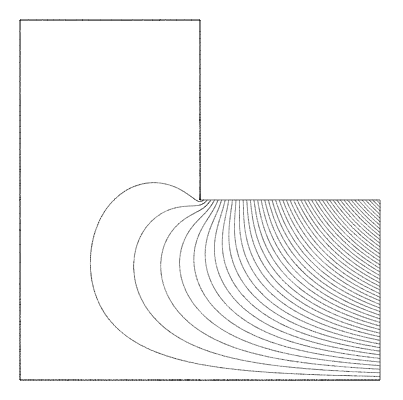} 
\includegraphics[width=1.8in]{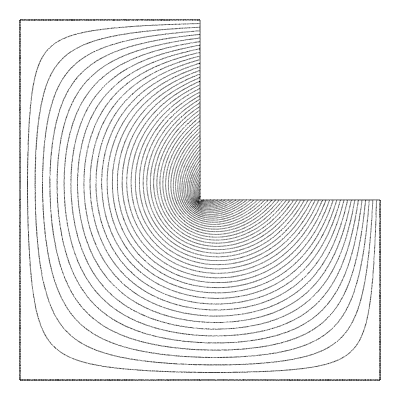} 
\end{center}
\caption{\label{LShapeElement} A ``degenerate'' octagon $K$, and contour
  plots of five of the eight basis functions for $V_1(K)$.}
\end{figure}

Let $\Omega\subset\RR^2$ be a polygon.
For a polygonal partition $\cT=\{K\}$ of $\Omega$ with 
vertices $\cV=\{z\}$ and edges
$\cE=\{e\}$, we use $\cV(K)$ and $\cE(K)$ to
denote, respectively, the vertices and edges of the polygon $K$.
Throughout, we use $\PP_j(S)$ to denote the polynomials of total degree at most $j$ on $S$, where $S$ is typically a polygon or a straight line segment, and we use the
convention that $\PP_j(S)=\{0\}$ when $j<0$.  We use $\PP_j(\partial K)$ to denote the continuous functions on
$\partial K$ which, when restricted to an edge $e\subset\partial K$,
are in $\PP_j(e)$.   We also briefly consider the space $\QQ_j(S)$ of polynomials of degree at most $j$ in each variable.
We allow degenerate polygons, i.e. those having a vertex (or more)
whose two adjacent edges form a straight angle, though we do not allow
two edges to meet at a zero angle (polygon with slit).  Allowing degenerate polygons
eliminates the possibility of ``hanging nodes'' in a polygonal
partition of $\Omega$.  A degenerate
octagon, congruent to an L-shaped hexagon, is shown in Figure~\ref{LShapeElement}. 
\begin{definition}[Shape-Regularity]\label{ShapeRegular}
A family of polygonal partitions $\cF=\{\cT\}$ is called shape-regular
when there are constants $c,\sigma>0$ such that, for every $\cT\in\cF$
and every $K\in\cT$:
\begin{enumerate}
\item $h_K\leq c h_e$ for all $e\in\cE(K)$, where
  $h_K=\mathrm{diam}(K)$ and $h_e=|e|$ is the length of the edge $e$.
\item $K$ is star-shaped with respect to a circle of (maximal) radius
  $\rho_K$, and $h_K\leq \sigma \rho_K$.
\end{enumerate}
Selecting such a circle of maximal radius, we may choose to denote its
center by $z_K$.
\end{definition}

\begin{definition}[Local Poisson Space]~\label{LocalSpace}
Given a polygon $K$ with $N$ edges/vertices, and an index $m\in\NN$, 
we define the local space $V_m(K)$ by
\begin{align}\label{VK}
  v\in V_m(K)\mbox{ if and only if } \Delta v\in\PP_{m-2}(K) \mbox{ in
  }K \mbox{ and } v\in\PP_m(\partial K)\mbox{ on }\partial K~.
\end{align}
It is clear that $\PP_m(K)\subset V_m(K)$, and $V_m(K)$ is naturally
decomposed as $V_m(K)=V^{K}_m(K)\oplus V^{\partial K}_m(K)$, where
\begin{align}\label{VKT}
  v\in V^{K}_m(K)\mbox{ if and only if } \Delta v\in\PP_{m-2}(K) \mbox{ in }K \mbox{ and } v=0\mbox{ on }\partial K~,\\
\label{VKE}
v\in V^{\partial K}_m(K)\mbox{ if and only if } \Delta v=0 \mbox{ in
}K \mbox{ and } v\in \PP_m(\partial K)\mbox{ on }\partial K~.
\end{align}
From this decomposition, it is apparent that
\begin{align}\label{dimVK}
\dim V_m(K)=\dim \PP_{m-2}(K) +\dim \PP_m(\partial K)=\binom{m}{2}+Nm~.
\end{align}
We may also decompose $V^{\partial K}_m(K)$ as
\begin{align}\label{VKSplit}
V^{\partial K}_m(K)=V_m^{\cV(K)}(K) \oplus  V_m^{\cE(K)}(K)~,
\end{align}
where $V_m^{\cV(K)}(K)=V_1(K)$, and $V_m^{\cE(K)}(K)$ consists of
those functions in $V_m^{\partial K}(K)$ that vanish at the vertices $\cV(K)$.
The decomposition $V_m(K) =V_m^{\cV(K)}(K) \oplus
V_m^{\cE(K)}(K)\oplus V_m^{K}(K) $ into \textit{ vertex, edge and
  interior functions}
corresponds naturally with the unisolvent set of degrees of freedom for $V_m(K)$,
\begin{align}\label{dofVK}
  v(z)\;\,\forall z\in\cV(K)\quad,\quad \int_e v p\,ds\;\,\forall p\in
  \PP_{m-2}(e)\;\,\forall e\in\cE(K) \quad,\quad \int_K
  vp\,dx\;\,\forall p\in \PP_{m-2}(K)~.
\end{align}
One might replace the moment-based edge degrees of freedom by
evaluations at $m-1$ distinct interior points on each edge, as
suggested, for example, in~\cite{vemMMMAS2013}.
\end{definition}

\begin{remark}\label{IntegralRelationsRemark} The following basic integral
  relations for the local Poisson space are often of use:
\begin{align}\label{IntegralRelations}
\int_K \nabla\varphi\cdot\nabla\phi\,dx=\begin{cases}
\int_{\partial K}\varphi\frac{\partial \phi}{\partial
  n}\,ds=\int_{\partial K}\phi\frac{\partial \varphi}{\partial
  n}\,ds&,\;\varphi,\phi\in V_m^{\partial K}(K)\\
-\int_K\varphi\Delta\phi\,dx=-\int_K\phi\Delta\varphi\,dx&,\, \varphi,\phi\in V_m^{K}(K)\\
0&,\; \varphi\in V_m^{\partial K}(K)\,,\, \phi\in V_m^{K}(K)
\end{cases}~.
\end{align}
For example, if the diffusion coefficient $A$ in~\eqref{ModelProblem}
is scalar and piecewise constant on $\cT$, the alternate forms of
$H^1$-inner-product above are typically employed in practice for the
formation of the finite element stiffness matrix.  We will also see in
Section~\ref{Interpolation}  how these integrals aid in the
understanding of interpolation in $V_m(K)$.
\end{remark}

\begin{remark}[Comparisons with $\PP_m(K)$ and
  $\QQ_m(K)$]\label{DimensionComparisons}
One sees that
\begin{align*}
\dim V_m(K)-\dim\PP_m(K)=(N-2)m-1\quad,\quad \dim Q_m(K)-\dim V_m(K)=m^2/2-(N-5/2)m+1~,
\end{align*}
and we already noted that $V_m(K)\supset \PP_m(K)$.
In the case of triangles ($N=3$), one immediately deduces that
$V_1(K)=\PP_1(K)$; but $\PP_m(K)$ is a proper subset of $V_m(K)$
when $N>3$ and/or $m>1$.  For quadrilaterals ($N=4$), $\dim
V_1(K)=\dim\QQ_1(K)$; but $V_1(K)\neq\QQ_1(K)$ for general
quadrilaterals, though they are the same for  rectangles aligning with the cardinal axes.
For generic $m$ and $N$, neither of these two spaces is contained in the other, and
$\QQ_m(K)$ typically has larger dimension.
\end{remark}

\begin{remark}[Singular Functions in $V_m(K)$]\label{SingularFunctions} 
  The typical singular behavior of functions in $V_m(K)$ near the
  corners of $K$ is well-understood
  (cf.~\cite{Grisvard1985,Grisvard1992,Wigley1964,Zargaryan1984}).
  Since $V_m(K)$ is finite dimensional, any basis we choose for this
  space must possess such singularities in some of its components.  For
  example, if we consider $V_m(K)$ for the degenerate octagon in
  Figure~\ref{LShapeElement}, and let $r=r(x)$ denote the distance
  from $x$ to the
  set of vertices, then the typical leading-singularity behavior of functions
 in this space is $r^2\ln r$ near each of the five vertices at
  $\pi/2$ internal angles, and $r^{2/3}$ near the vertex
  at the $3\pi/2$ internal angle.
\end{remark}

\begin{definition}[Global Poisson Space]\label{GlobalSpace}
The global Poisson space corresponding to the partition $\cT=\{K\}$ is
\begin{align}\label{VSpace}
V_m=V_m(\cT)=\{v\in C(\overline\Omega):\,v\vert_K\in V_m(K)\mbox{ for each } K\in
\cT\}~,
\end{align}
and a unisolvent set of degrees of freedom is given by
\begin{align}\label{dofV}
v(z)\;\,\forall z\in\cV\quad,\quad \int_e v p\,ds\;\,\forall p\in \PP_{m-2}(e)\;\,\forall e\in\cE
\quad,\quad \int_K vp\,dx\;\,\forall p\in \PP_{m-2}(K)\;\,\forall K\in\cT~.
\end{align}
As with the local Poisson spaces, we may naturally decompose $V_m$
into \textit{vertex, edge and interior functions},
\begin{align}\label{VmDecomposition}
V_m=V_m^\cV\oplus V_m^\cE\oplus V_m^\cT~.
\end{align}
We clearly have
\begin{align}\label{dimV}
\dim V_m={\rm card}(\cT)\,\binom{m}{2}+{\rm card}(\cE)\,(m-1)+{\rm card}(\cV)~,
\end{align} 
and dimension of the
space $V_m\cap\cH$ is suitably reduced by replacing $\cV$ and $\cE$ in~\eqref{dimV},
with the non-Dirichlet vertices $\cV'\not\subset \partial\Omega_D$ and
edges $\cE'\not\subset\partial\Omega_D$. 
\end{definition}

\begin{remark}\label{MixedSpace}
  An obvious variant of the global space described above might include
  a mixture standard finite elements on triangles and/or rectangles
  throughout much of the domain, connected to the polygonal Poisson
  elements by matching polynomial basis functions along shared edges.
  The restrictions that a computational cell $K$ is a polygon and
  that the boundary data is piecewise polynomial on $\partial K$ may
  be relaxed as well, again provided that there is a convenient mechanism
  enforcing agreement at interfaces between elements.  As will be seen
  below, the Nystr\"om approach readily provides the flexibility to
  explore such variants.
\end{remark}

\subsection{Integral representations of functions in $V_m(K)$}

\begin{remark}[Polynomial Solutions of Poisson Problems with
  Polynomial Sources]\label{KarachikAntropova}
  Suppose that $p\in\PP_j(\RR^n)$ is homogeneous and of degree $j$,
  i.e. $p(c x)=c^jp(x)$, and define $q\in\PP_{j+2}(\RR^n)$ by
\begin{align}\label{KarachikAntropovaEq}
 q(x)&=\sum_{k=0}^{[j/2]}
\frac{(-1)^k \Gamma(n/2+j-k)}{\Gamma(n/2+j+1)(k+1)!}\left(\frac{|x|^2}{4}\right)^{k+1}\,\Delta^k p(x)~,
\end{align}
where $[j/2]$ denotes the integer part of $j/2$.  It is shown
in~\cite[Theorem 2]{KarachikAntropova2010} that
$\Delta q = p$. Now recall that $v\in V_m^K(K)$ satisfies $\Delta v=p$
in $K$ for some $p\in \PP_{m-2}(K)$, with
$v=0$ on $\partial K$.  So we see that there is a $q\in \PP_{m}(K)$
for which $w=v-q$ satisfies $\Delta w=0$ in $K$, with
$w=-q\in\PP_m(\partial K)$ on $\partial K$.
\end{remark}

A practical consequence of Remark~\ref{KarachikAntropova} is that the computation of any function in
$V_m(K)=V^K_m(K)\oplus V^{\partial K}_m(K)$
is reduced to solving problems of the form
\begin{align}\label{LaplaceEquation}
\Delta w=0\mbox{ in }K\quad,\quad w=g\in\PP_m(\partial K)\mbox{ on }\partial K~.
\end{align}
Let $G(x,y)$ denote the fundamental solution for the Laplacian,
\begin{align}\label{FundamentalSolution}
G(x,y)=\frac{1}{2\pi}\,\ln\frac{1}{|x-y|}~.
\end{align}
We may express the solution of~\eqref{LaplaceEquation} as a double-layer potential,
\begin{align}\label{DoubleLayer}
w(x)=\int_{\partial K} \frac{\partial G(x,y)}{\partial n(y)}\,\phi(y)\,ds(y)
\mbox{ for } x\in K~,
\end{align}
where the density $\phi$ satisfies the second-kind integral equation
\begin{align}\label{SecondKind}
\frac{\phi(x)}{2}-\int_{\partial K} \frac{\partial G(x,y)}{\partial n(y)}\,\phi(y)\,ds(y)=-g(x)\mbox{ for }x\in\partial K~.
\end{align}
Here and following, $n=n(y)$ denotes the outward unit normal at $y$ to the
domain under consideration.  The basic approach of this paper is to
approximate $w\in V_m(K)$ via~\eqref{DoubleLayer}-\eqref{SecondKind}
by Nystr\"om discretizations, as discussed in Section~\ref{Nystrom}.
In contrast, the BEM-FEM approach expresses $w$ as a combination of single- and double-layer potentials,
\begin{align}\label{SingleLayer}
w(x)=\int_{\partial K} G(x,y)\,\psi(y)\,ds(y)-\int_{\partial K} \frac{\partial G(x,y)}{\partial n(y)}\,g(y)\,ds(y)
\mbox{ for } x\in K~,
\end{align}
where the density $\psi$ satisfies the first-kind integral equation
\begin{align}\label{FirstKind}
\int_{\partial K} G(x,y)\,\psi(y)\,ds(y)=\frac{g(x)}{2}+\int_{\partial K} \frac{\partial G(x,y)}{\partial n(y)}\,g(y)\,ds(y)
\mbox{ for } x\in \partial K~.
\end{align}

If~\eqref{SecondKind} and~\eqref{FirstKind} are to be
understood pointwise, then the left-hand side of~\eqref{SecondKind}
and the right-hand side of~\eqref{FirstKind} must be modified
at each corner $z\in\partial K$ (cf.~\cite{KressBook2014}).  For
example, if the interior angle at $z$
is $\alpha\pi$,
then~\eqref{SecondKind} becomes
 \begin{align}\label{SecondKindCornerA}
\left(1-\frac{\alpha}{2}\right)\phi(z)-\int_{\partial K} \frac{\partial G(z,y)}{\partial n(y)}\,\phi(y)\,ds(y)=-g(z)~.
\end{align}
In practice, it is more convenient to use a modified form that
provides a smoother integrand and does not require specific knowledge
of the angle at $z$.  In particular, if we take $\partial K'$ to
denote $\partial K$ without the corners, we have
 \begin{align}\label{SecondKindCornerB}
\frac{\phi(x)+\phi(z)}{2}-\int_{\partial K} \frac{\partial
  G(x,y)}{\partial n(y)}(\phi(y)-\phi(z))\,ds(y)=-g(x)\mbox{ for
}x\in\partial K' \cup\{z\}~.
\end{align}

\begin{remark}\label{Densities}
The derivation of~\eqref{SingleLayer}-\eqref{FirstKind} from Green's
formulas reveals that $\psi=\partial w/\partial n$, so the solution of
\eqref{FirstKind} directly provides a \textit{Dirichlet-to-Neumann
  Map} $g\mapsto \partial w/\partial n$.  In contrast, the density
$\phi$ in~\eqref{SecondKind} is given by $\phi=v-g$, where $v$ is the
unique solution of the complementary exterior Neumann problem
\begin{align*}
\Delta v=0\mbox{ in }\RR^2\setminus\overline K\quad,\quad
\frac{\partial v}{\partial n}=\frac{\partial w}{\partial n}~,
\end{align*}
and $v=o(1)$ uniformly in all directions as $|x|\to\infty$.
As with the discussion of singular behavior in $V_m(K)$ in Remark~\ref{SingularFunctions},
the singular behavior of the densities $\phi,\psi$ near corners is
well-understood.  Let $z\in\partial K$ be a corner with interior angle $\alpha\pi$,
and let $x\in\partial K$ be near $z$.  If
$\alpha\in(0,2)$ is irrational, we have
\begin{enumerate}
\item $\psi(x)\sim |x-z|^{\alpha^{-1}-1}$ as $|x-z|\to 0$.  So $\psi$
  is expected to blow up near non-convex corners ($\alpha\in(1,2)$),
  and its tangential derivative is expected to blow up at convex
  corners ($\alpha\in(0,1)$).
\item $\phi(x)-\phi(z)\sim |x-z|^{\sigma^{-1}}$ as $|x-z|\to 0$, where
  $\sigma=\max(\alpha,2-\alpha)\geq 1$.  So $\phi$ will be bounded at
  all corners, but its tangential derivative will typically blow up at
  each non-straight corner.  In the case that $w$ is smooth in $K$,
  the complementary $v$ does not ``inherit'' any singular behavior
  from $w$ via the Neumann boundary condition, and we have $\sigma=2-\alpha$. 
\end{enumerate}
In the case of rational $\alpha$, logarithmic terms may appear in the
asymptotic expansions of $\phi$ and $\psi$, but they are only the
dominant terms in the expansion when: $\alpha =1/2$ for $\psi$, where
we have $\psi(x)\sim |x-z|\,\ln |x-z|$ as $|x-z|\to 0$; or $\alpha =3/2$ and $w$
is smooth, where $\phi(x)-\phi(z)\sim |x-z| ^2\,\ln |x-z|$ as $r\to 0$.  In
any case, $\phi$ is H\"older continuous.
\end{remark}

\subsection{A Hierarchical Basis for $V_m(K)$}
The decomposition $V_m(K)=V_m^K(K)\oplus V_m^{\partial K}(K)$ makes it
clear that a basis for $\PP_{m-2}(K)$ yields a corresponding basis for
$V_m^K(K)$, and a basis for $\PP_{m}(\partial K)$ yields a corresponding basis for
$ V_m^{\partial K}(K)$.  In light of Remark~\ref{KarachikAntropova}, it is convenient to choose
a basis for $\PP_{m-2}(K)$ in terms of translated monomials, centered
at some convenient point $z\in K$,  
\begin{align}\label{MonomialBasis}
\PP_{m-2}(K)=\mbox{span}\{(x-z)^\beta=(x_1-z_1)^{\beta_1}
(x_2-z_2)^{\beta_2}:\;|\beta|=\beta_1+\beta_2\leq m-2\}~.
\end{align}
Such bases are naturally hierarchical in polynomial degree.  

Given an edge $e\subset \partial K$, with endpoints $z,z'$, we construct a hierarchical basis
for $\PP_m(e)$ as follows.  Let
$\lambda_z^e,\lambda_{z'}^e\in\PP_1(e)$ be the corresponding
barycentric coordinates for $e$, defined by
$\lambda_z^e(z')=\lambda_{z'}^e(z)=\delta_{zz'}$.  For $2\leq j\leq
m$, we define $b_j^e(x)=\hat{L}_j(\lambda_z^e(x)-\lambda_{z'}^e(x))$,
where $\hat{L}_j$ is the integrated Legendre polynomial of degree $j$
(cf.~\cite{BabushkaSzaboBook}).  These are given in terms of the
standard Legendre polynomials $L_i$, with normalization $L_i(1)=1$, by
$\hat{L}_j(t)=\int_{-1}^tL_{j-1}(s)\,ds=(L_j(t)-L_{j-2}(t))/(2j-1)$.
We see that $b_j^e(z)=b_j^e(z')=0$, and a basis for $\PP_m(e)$ that is
hierarchical in polynomial degree is
\begin{align}\label{HierarchicalEdgeBasis}
\{\lambda_z^e,\lambda_{z'}^e\}\cup\{b_2^e,\ldots,b_{m}^e\}\doteq
\cB_1^e\cup\cB_m^e~.
\end{align}
Given a vertex $z\in\cV(K)$, the function $\lambda_z\in\PP_1(\partial
K)$ is determined by the conditions $\lambda_z(z')=\delta_{zz'}$ for all
$z'\in\cV(K)$; so $\lambda_z\vert_e=\lambda_z^e$ if $z$ is an endpoint
of $e$. Clearly $\cB_1=\{\lambda_z:\;z\in\cV(K)\}$ is a basis for
$\PP_1(\partial K)$.  Each element $b_j^e$ of $\cB_m^e$ vanishes at the
endpoints of $e$, so we continuously extend it by $0$ to $\partial K$.
Finally, a hierarchical basis for
$\PP_m(\partial K)$, is given by
\begin{align}
\cB_1\cup\left(\bigcup_{e\in\cE(K)}\cB_m^e\right)\doteq \cB_1\cup\cB_m~.
\end{align}

\section{Nystr\"om Approximation Second-Kind Integral Equations}\label{Nystrom}

As was seen in the previous section, the computation of $v\in V_m(K)$
is reduced to the computation of a harmonic function $w$ on $K$ with
prescribed Dirichlet data, and we opt to do so via second-kind
integral equations.  Nystr\"om methods~\cite{Nystrom1928,Nystrom1930}
for second-kind integral equations, in their most basic forms, are
derived by replacing the boundary integral with a suitable quadrature,
and sampling the resulting equation at the quadrature points.  The
performance of the method is directly tied to the performance of the
underlying quadrature, and we
will briefly describe the version proposed by Kress~\cite{Kress1990} for
problems of the sort that we here consider, after first looking more
closely at the components of the integrand.

Recalling~\eqref{SecondKindCornerB},
for $x\in\partial K$ near or at a vertex/corner $z$, we have
\begin{align}\label{NearVertex}
\frac{\phi(x)+\phi(z)}{2}+\int_{\partial K}
F(x,y)(\phi(y)-\phi(z))\,ds(y)=-g(x)\;,\;
F(x,y)=-\frac{\partial G(x,y)}{\partial n(y)}=-\frac{(x-y)\cdot n(y)}{2\pi\,|x-y|^2}~.
\end{align}
We note that $F(x,y)=0$ when $x$ and $y$ are on the same (straight) edge of $\partial K$.
More generally, for any fixed $x\in\partial K$, $F(x,y)$ is a piecewise smooth function of $y\in\partial K$, with
bounded jump-discontinuities at the corners of $\partial K$.  In fact,
for any $x\in\partial K$, 
$F(x,y)$ is analytic
in the interior of each edge.  This is not to say that $F(x,y)$ exhibits
no difficult behavior: if $x\in\partial K$ is very near (but not at) a
corner $z$, then $F(x,y)$ and its tangential derivatives in $y$ are very large as
$y$ approaches $z$ along the edge not containing $x$.  More
specifically, if $x$ and $y$ are on opposite straight edges sharing $z$, and
the interior angle at $z$ is $\alpha\pi$
then
\begin{align*}
\lim_{y\to z} F(x,y)=\frac{\sin(\alpha\pi)}{2\pi|x-z|}~,
\end{align*}
where $y$ is understood to approach $z$ along the edge they share.
For $x\in\partial K$ we choose $z=z(x)$ to be the nearest vertex in
terms of distance along the boundary, breaking ties arbitrarily if $x$
is at the midpoint of an edge.  The fact that the integrand in~\eqref{NearVertex} vanishes at $y=z$
makes it easier to approximate the integral by simple quadrature.

The basic quadrature employed by Kress~\cite{Kress1990}  for $f\in C[0,1]$ is obtained by applying the uniform
trapezoid rule after a sigmoidal change-of-variable,
\begin{align*}
\int_0^1f(t)\,dt=\int_0^1f(\eta(\tau))\eta'(\tau)\,d\tau\approx \frac{1}{n}\sideset{}{''}\sum_{k=0}^nf(\eta(k/n)) \eta'(k/n)=\frac{1}{n}\sum_{k=1}^{n-1}f(\eta(k/n)) \eta'(k/n)~,
\end{align*}
where the transformation $t=\eta(\tau)$ is given by
\begin{align*}
\eta(\tau)=\frac{[c(\tau)]^p}{[c(\tau)]^p+[1-c(\tau)]^p} \quad,\quad c(\tau)=\left(\frac{1}{2}-\frac{1}{p}\right)
(2\tau-1)^3+\frac{1}{p}(2\tau-1)
+\frac{1}{2}~,
\end{align*}
and $p\geq 2$ is an integer.  It is straight-forward to see that $\eta$ has a root of order $p$ at $0$, and
$1-\eta$ has a root of order $p$ at $1$.  A careful convergence analysis of this quadrature is given in~\cite{Kress1990}, 
showing that it is convergent on $C[0,1]$, and, for the kinds of integrands
we encounter here, of increasingly higher-order in $n$ as $p$ is increased.
If $e$ is a smooth (curved) edge, with smooth 
parametrization $x_e=x_e(t)$ satisfying $|x_e'(t)|\geq \sigma>0$, we have
the quadrature 
\begin{align}\label{KressQuad}
\int_e f\,ds=\int_0^1 f(x(t))|x_e'(t)|\,dt\approx \sideset{}{''}\sum_{k=0}^nf(x_k^e)
\omega_k=\sum_{k=1}^{n-1}f(x_k^e) \omega_k~,
\end{align}
where  $\omega_k=\eta'(k/n) \, |x_e'(\eta(k/n))|/n$ and $x_k^e=x_e(\eta(k/n))$.
In the case of a straight edge $e$ having endpoints $z,z'$,
these weights and points simplify to 
\begin{align}\label{KressQuadStraight}
\omega_k= \eta'(k/n)\,|e|/n\quad,\quad x_k^e=(1-\eta(k/n)) z+\eta(k/n)z'=\eta(1-k/n) z+\eta(k/n)z'~.
\end{align}

Keeping a fixed $n$ and $p$ for all edges, we take a global enumeration of the quadrature points and weights
(including vertices), $\{(x_j,\omega_j):\,1\leq j\leq M=nN\}$.  The
Nystr\"om linear system corresponding
to~\eqref{NearVertex} is given by
\begin{align}
\frac{\phi_i+\phi_k}{2}+\sum_{j=1}^M
F(x_i,x_j)(\phi_j-\phi_k)\omega_j=-g(x_i)~,\label{NystromSystem}
\end{align}
where $x_k=z(x_i)$ is the vertex nearest $x_i$.
The approximation $\tilde{w}(x)\approx w(x)$ for
$x\in K$ is
given by
\begin{align}\label{wTilde}
\tilde{w}(x)=-\sum_{j=1}^M F(x,x_j)\phi_j\omega_j~.
\end{align}

We demonstrate the efficacy of the Nystr\"{o}m scheme in obtaining
accurate approximations to solutions of a couple of example boundary
value problems~\eqref{LaplaceEquation}
that present challenges similar to those that arise in
the construction of $V_m(K)$.  In particular, these examples deal with
regions that have multiple corners, so the corresponding densities $\phi$ exhibit
singularities as described in Remark~\ref{Densities}.  In both cases
a harmonic function $w$ is given, and the (non-polynomial) Dirichlet
data is taken from $w$.  Relative and/or
absolute errors in the Nystr\"{o}m
approximation of $w$ are given at several points in the interior of
$K$ for these examples, when $n$ points are used on each edge and the
parameter $p=6$ is used for the quadrature. 

\begin{example}
Let $K$ be the L-shaped hexagon with vertices at $(0,0), (1,0), (1,1), (-1,1),
(-1,-1), (0,-1)$, and take $w=\ln|x-\hat{x}|$, where $\hat{x}=(10,0)$.  Although $w$ is smooth in $K$, the corresponding density $\phi$ will
have singular behavior as discussed before.  More specifically,
$\psi\sim |x-z|^{2/3}$ near each of the five corners $z$ having
interior angle $\pi/2$, and $\psi\sim |x|^{2}\ln|x|$ near the corner
at the origin having interior angle $3\pi/2$.  In
Table~\ref{table::conv::L}, relative errors in the Nystr\"{o}m
approximation of $w$ are given at four points in $K$, and
clearly demonstrate the high-order convergence as the number of points
per edge increases, as well as the accuracy even when few points are
used.

\begin{table}
\caption{Relative errors at five points for the Nystr\"om approximation of
  $w=\ln|x-\hat{x}|$ in an L-shaped hexagon.}
\label{table::conv::L}
\begin{tabular}{lccccc}
\hline
$n$ & { (0.5,0.5) } & { (0.1,0.1)} & { (0.01,0.01)} & { (0.001,0.001)} & { (0.999,0.001)} \\\hline
16 & 5.954e-07 & 1.168e-05 & 3.231e-06 & 3.142e-07 &1.912e-05\\
32 & 1.077e-10 & 1.976e-07 & 2.530e-08 & 3.379e-07 &1.298e-06\\
64 & 6.565e-13 & 1.628e-09 & 5.423e-10 & 3.584e-09 &1.867e-08\\
128 &  9.857e-15 & 4.329e-11 & 2.062e-11 & 1.343e-11&2.897e-10 \\
256 & 1.971e-16 & 4.990e-13 & 9.927e-13 & 1.341e-12 &1.309e-11\\
512 & 0.000e+00 & 3.680e-15 & 7.120e-14 & 7.175e-14 &6.532e-13\\
1024 & 0.000e+00 & 5.811e-16 & 4.631e-15 & 7.522e-15 &3.941e-14\\
2048 & 0.000e+00 & 0.000e+00 & 1.929e-16 & 1.929e-16 &5.457e-15\\
\hline
\end{tabular}
\end{table}
\end{example}

\begin{example}\label{NystromPacman}
Here we take $K=K(\alpha)= \{ x=(r\cos\theta, r\sin\theta) \, :\;\ 0 < r
< 1, 0 < \theta < \alpha\pi \}$ to be the sector of the unit circle
with interior opening angle $\alpha \pi$, and
$w=r^{1/\alpha}\sin(\theta/\alpha)$.  This example provides a
situation like many we expect to encounter in practice, where both
$\phi$ and $w$ have singular behavior near a corner.
The interior angles at the
other two corners, $(0,1)$ and $(\cos \alpha \pi, \sin \alpha \pi)$, are both
$\pi/2$.  We consider the case
$K(3/2)$, a ``circular L-shape''.  We have $\phi(x)\sim |x-z|^{2/3}$ near each
of the three corners $z$.  Since $w=0$ at the origin, we provide both 
relative and absolute approximation errors at a few points near the origin.
The results for 
$K(3/2)$ are given in Table \ref{table::conv::pacman}.  We again see
similar high-order convergence.

\begin{table}
  \caption{Relative and absolute errors at three points for the
    Nystr\"om approximation of $w=r^{2/3}\sin(2\theta/3)$ on the
    circular L-shape $K(3/2)$.}
  \label{table::conv::pacman}
\begin{tabular}{lcccccc}
\hline
& \multicolumn{2}{c}{(0.1,0.1)}& \multicolumn{2}{c}{(0.01,0.01)} & \multicolumn{2}{c}{(0.001,0.001)}\\
$n$ & Rel & Abs &  Rel & Abs & Rel & Abs \\\hline
16 & 9.518e-03 & 1.292e-03 & 3.065e-02 & 8.963e-04 & 5.558e-02 & 3.501e-04 \\
32  & 2.841e-05 & 3.856e-06 & 1.237e-02 & 3.616e-04 & 8.311e-02 & 5.236e-04 \\
64  & 2.763e-08 & 3.750e-09 & 2.910e-05 & 8.510e-07 & 3.303e-05 & 2.228e-07 \\
128 & 9.083e-11 & 1.233e-11 & 1.715e-09 & 5.016e-11 &  2.008e-06 & 1.265e-08 \\
256 & 2.708e-13 & 3.675e-14 & 7.337e-12 & 2.145e-13 & 1.590e-10 & 1.002e-12 \\
512 & 8.180e-16 & 1.110e-16 & 4.675e-14 & 1.367e-15 & 8.468e-13 & 5.334e-15 \\
1024 &  3.190e-14 & 4.330e-15 & 2.442e-13 & 7.140e-15 & 9.173e-13 & 5.778e-15 \\
2048 & 2.699e-14 & 3.664e-15 & 2.442e-13 & 7.140e-15 & 1.339e-12 & 8.432e-15 \\
\hline
\end{tabular}
\end{table}
\end{example}

\section{Interpolation in $V_m$}\label{Interpolation}

We consider some properties of the interpolation operator
$\cI_m:C(\overline\Omega)\to V_m$ defined by
\begin{align}\label{InterpolationOperator}
\begin{cases}\left(\cI_m v\right)(z)=v(z)&\;\,\forall z\in\cV\\
\int_e \left(\cI_m v\right) p\,ds=\int_e v p\,ds&\;\,\forall p\in \PP_{m-2}(e)\;\,\forall e\in\cE\\
 \int_K \left(\cI_m v\right)p\,dx=\int_K v p\,dx&\;\,\forall p\in \PP_{m-2}(K)\;\,\forall K\in\cT
\end{cases}~,
\end{align}
both theoretically and empirically.
The interpolation operator $\cI_m$ may be decomposed in such a way as
to correspond to the space decomposition $V_m= V_m^\cV\oplus
V_m^\cE\oplus  V_m^\cT$, namely
$\cI_m=\cI_m^\cV+\cI_m^\cE+\cI_m^\cT$, where
$\cI_m^\cV:C(\overline\Omega)\to V_m^\cV$, 
$\cI_m^\cE:C(\overline\Omega)\to V_m^\cE$ and $
\cI_m^\cT:C(\overline\Omega)\to V_m^\cT$ are uniquely determined by
\begin{align}\label{InterpolationOperatorB}
\begin{cases}\left(\cI_m^\cV v\right)(z)=v(z)&\;\,\forall z\in\cV\\ \int_e \left(\cI_m^{\cE}v\right) p\,ds=\int_e (v-\cI_m^{\cV}v) p\,ds&\;\,\forall p\in \PP_{m-2}(e)\;\,\forall e\in\cE\\
 \int_K \left(\cI_m^\cT v\right)p\,dx=\int_K \left(v-\cI_m^\cV v-\cI_m^\cE v\right)p\,dx&\;\,\forall p\in \PP_{m-2}(K)\;\,\forall K\in\cT
\end{cases}~.
\end{align}
We have the obvious restrictions of these interpolation operators to a single
element, and we use the same symbols to denote them.  On a single
element  $K$, it is also convenient to use $\cI_m^{\partial
  K}:C(\overline{K})\to V_m^{\partial K}(K)$ to denote the
interpolation operator defined by
\begin{align}\label{InterpolationOperatorC}
\begin{cases}\left(\cI_m^{\partial K} v\right)(z)=v(z)&\;\,\forall
  z\in\cV(K)\\ \int_e \left(\cI_m^{\partial K}v\right) p\,ds=\int_e (v-\cI_m^{\cV}v) p\,ds&\;\,\forall p\in \PP_{m-2}(e)\;\,\forall e\in\cE(K)
\end{cases}~.
\end{align}
So $\cI_m=\cI_m^{K}+\cI_m^{\partial K}$ on $C(\overline{K})$.

Given $v\in C(\overline{K})$, with $v=g$ on $\partial K$, we define
$\tilde{g}\in C(\partial\Omega)$ by
\begin{align}\label{InterpolationOnEdges}
\tilde{g}(z)=g(z)\;\forall z\in\cV(K)\quad,\quad \int_e \tilde{g}
p\,ds=\int_e g p\,ds\;\forall p\in\PP_{m-2}(e)\;\forall e\in\cE(K)~.
\end{align}
By definition, $\cI_m^{\partial K}v$ is the solution of 
\begin{align}\label{HarmonicInterpolant}
\Delta(\cI_m^{\partial K}v)=0\mbox{ in }K\quad,\quad \cI_m^{\partial
  K}v=\tilde{g}\mbox{ on }\partial K~.
\end{align}
Assuming that $g$ is continuously differentiable along each edge of $K$, 
we see that, for any $q\in\PP_m(e)$,
\begin{align*}
\int_e \frac{\partial(g-\tilde{g})}{\partial t}\frac{\partial
  q}{\partial t}\,ds=\left.(g-\tilde{g}) \frac{\partial q}{\partial
    t}\right\vert_{z'}^z-\int_e(g-\tilde{g})\frac{\partial^2 q}{\partial t^2}\,ds=0~,
\end{align*}
where $z,z'$ are the endpoints of $e$, and the partial derivatives are
in the tangential direction.  From this orthogonality relation, it is
clear that
\begin{align*}
|g-\tilde{g}|_{H^1(e)}=\inf_{p\in\PP_m(e)} |g-p|_{H^1(e)}\;\forall
e\in\cE(K)~.
\end{align*}
Recalling our hierarchical basis~\eqref{HierarchicalEdgeBasis} for $\PP_m(e)$, we have
\begin{align}\label{tildeg}
\tilde{g}(x)&=g(z)\lambda_z^e(x)+g(z')\lambda_{z'}^e(x)+\sum_{j=2}^mc_j
b_j^e(x)\quad,\quad
c_j=\frac{(2j-1)}{2}\int_e \frac{\partial g}{\partial t}\, L_{j-1}(\lambda_z-\lambda_{z'})\,ds~,
\end{align}
on $e$.  This follows from the relations
\begin{align*}
\frac{\partial b_j^e}{\partial
  t}=\frac{2}{|e|}\,L_{j-1}(\lambda_z-\lambda_{z'})\quad,\quad 
\int_e \frac{\partial b_j^e}{\partial  t}\frac{\partial b_i^e}{\partial  t}\,ds=\frac{4\delta_{ij}}{(2j-1)|e|}~.
\end{align*}

Having computed $\cI_m^{\partial K}u$
from~\eqref{HarmonicInterpolant}, with $\tilde{g}$ given
by~\eqref{tildeg}, we obtain $\cI_m v =\cI_m^{K}v+\cI_m^{\partial K}v$
on $K$ via 
\begin{align}\label{InteriorInterpolant}
\int_K \left(\cI_m^K v\right)p\,dx=\int_K \left(v-\cI_m^{\partial K} v\right)p\,dx&\;\,\forall p\in \PP_{m-2}(K)~.
\end{align}
Finally, taking a translated monomial basis $\{p_\beta\in \PP_{m-2}:\,|\beta|\leq
m-2\}$ of $\PP_{m-2}$, as in~\eqref{MonomialBasis}, and letting
$\{\varphi_\beta\in H^1_0(\Omega):\,\Delta\varphi_\beta=p_\beta\;,\; |\beta|\leq
m-2\}$ be the associated basis of $V_m^K(K)$, we see that the
coefficients $c_\beta$ of $\cI_m^K v=\sum_{|\beta|\leq 2}c_\beta
\varphi_\beta$ satisfy the negative definite linear system
\begin{align}\label{InteriorInterpolantSystem}
\sum_{|\beta|\leq m-2}c_\beta\int_K \varphi_\beta p_{\beta'}\,dx=\int_K
  \left(v-\cI_m^{\partial K} v\right) p_{\beta'}\,dx\;,\;\forall
  |\beta'|\leq m-2~.
\end{align}
Using~\eqref{IntegralRelations}, this system matrix is seen to be negative definite by 
\begin{align*}
\int_K \varphi_\beta p_{\beta'}\,dx=-\int_K \nabla\varphi_\beta \cdot\nabla\varphi_{\beta'}\,dx~.
\end{align*}

\begin{remark}\label{HarmonicApprox}
 When $m=1$, $\cI_m v =\cI_m^{\partial K}v$, so
interpolation in this case reduces to solving harmonic problems with
interpolated boundary data~\eqref{tildeg}.  If $v$ is harmonic in
$K$, then it can be well-approximated by functions in
$V_m^{\partial K}(K)$, because this space contains the harmonic
polynomials of degree $\leq m$, and these approximate general harmonic
functions essentially as well as the entire space $\PP_m(K)$ does (cf.~\cite{Melenk1997,Melenk1999}).  
\end{remark}

\begin{remark}\label{SolvingForGTilde}
  Although~\eqref{tildeg} provides a exact expression for the
  coefficients $c_j$ of $\tilde{g}$ on $e$ as an integral involving
  its tangential derivative, we found it more convenient for our interpolation
  experiments below to compute these coefficients in a different way.
  They were computed as the solution of a simple linear system derived
  by plugging the expression for $\tilde{g}$ from~\eqref{tildeg} into
  the integral identities in~\eqref{InterpolationOnEdges}, with $p$
  being the Lagrange polynomials of degree $\leq m-2$ associated with
  the edge.  The orthogonality relations for Lagrange polynomials
  leads to an $(m-1)\times(m-1)$ system matrix whose only non-zero
  elements are on the main diagonal and second lower-diagonal, and the
  righthand side is adequately addressed by quadrature.

\end{remark}

\begin{example}\label{Example1}
  To demonstrate the interpolation properties, we consider a numerical
  experiment and interpolate the function $v(x)=\sin(2\pi
  x_1)\sin(2\pi x_2)$ over $\Omega=(0,1)^2$ on a sequence of uniformly
  refined meshes, see Figure~\ref{fig:INT_Meshes}. The expansion
  coefficients are determined as described above. The volume integrals
  in~\eqref{InteriorInterpolantSystem} are realized by means of
  numerical quadrature. For this purpose, we split the element~$K$
  into triangles by connecting the vertices~$\cV(K)$ with the center
  of mass. Afterwards, a $7$-point Gaussian rule is applied on each
  triangle and the discrete functions in $V_m(K)$ are treated by means
  of Nystr\"om approximations, see Section~\ref{Nystrom}. The relative
  interpolation error is plotted in Figure~\ref{fig:INT_Convergence}
  for the $L_2$- as well as the $H^1$-norm with respect to the maximal
  mesh size $h_\mathrm{max} = \max\{h_K:K\in\cT\}$. An optimal order
  of convergence is achieved.
\begin{figure}[tbp]
 \centering
 \includegraphics[trim=4.4cm 7.4cm 4.4cm 7.8cm, width=0.22\textwidth, clip]{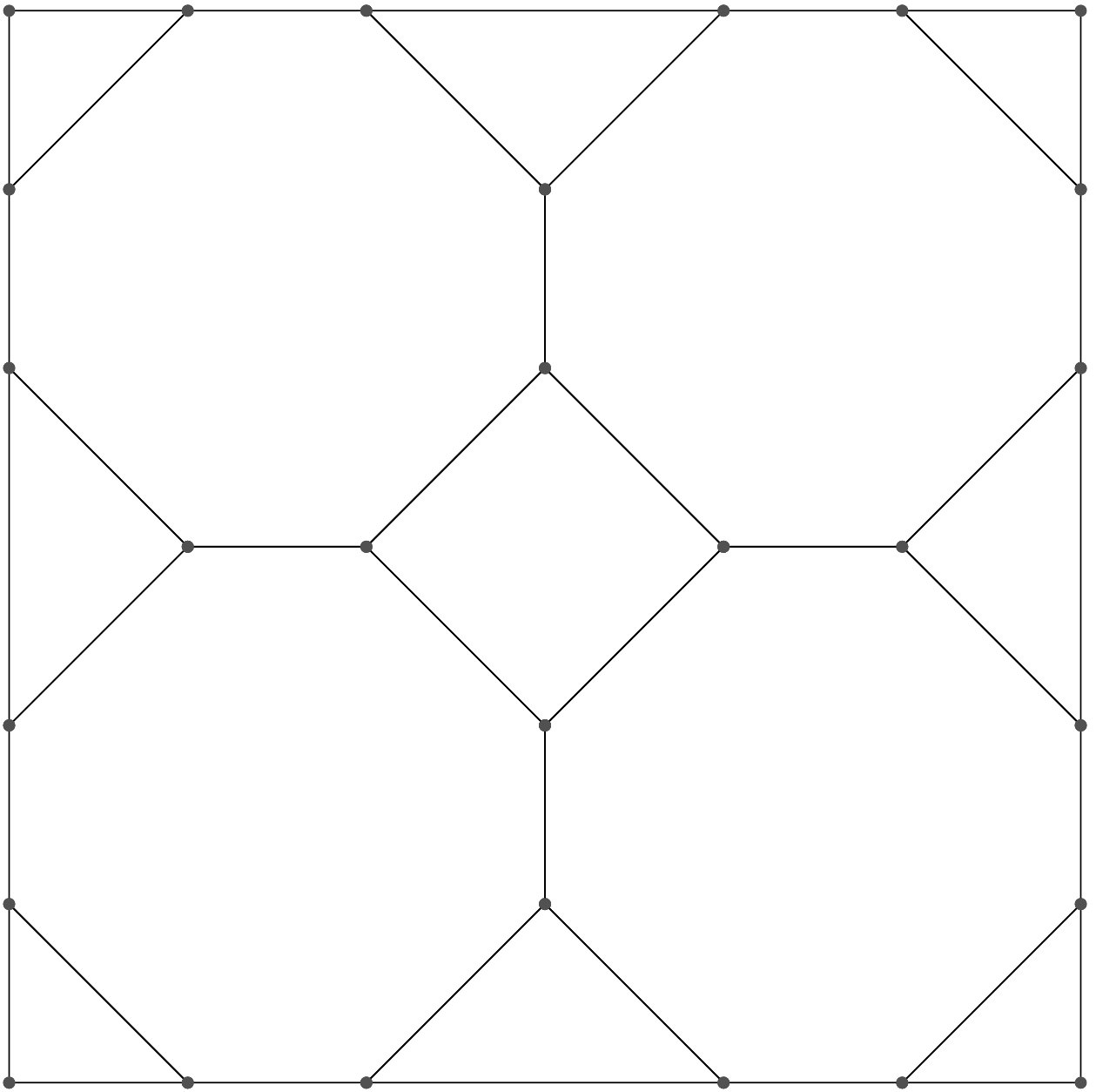}\hfill
 \includegraphics[trim=4.4cm 7.4cm 4.4cm 7.8cm, width=0.22\textwidth, clip]{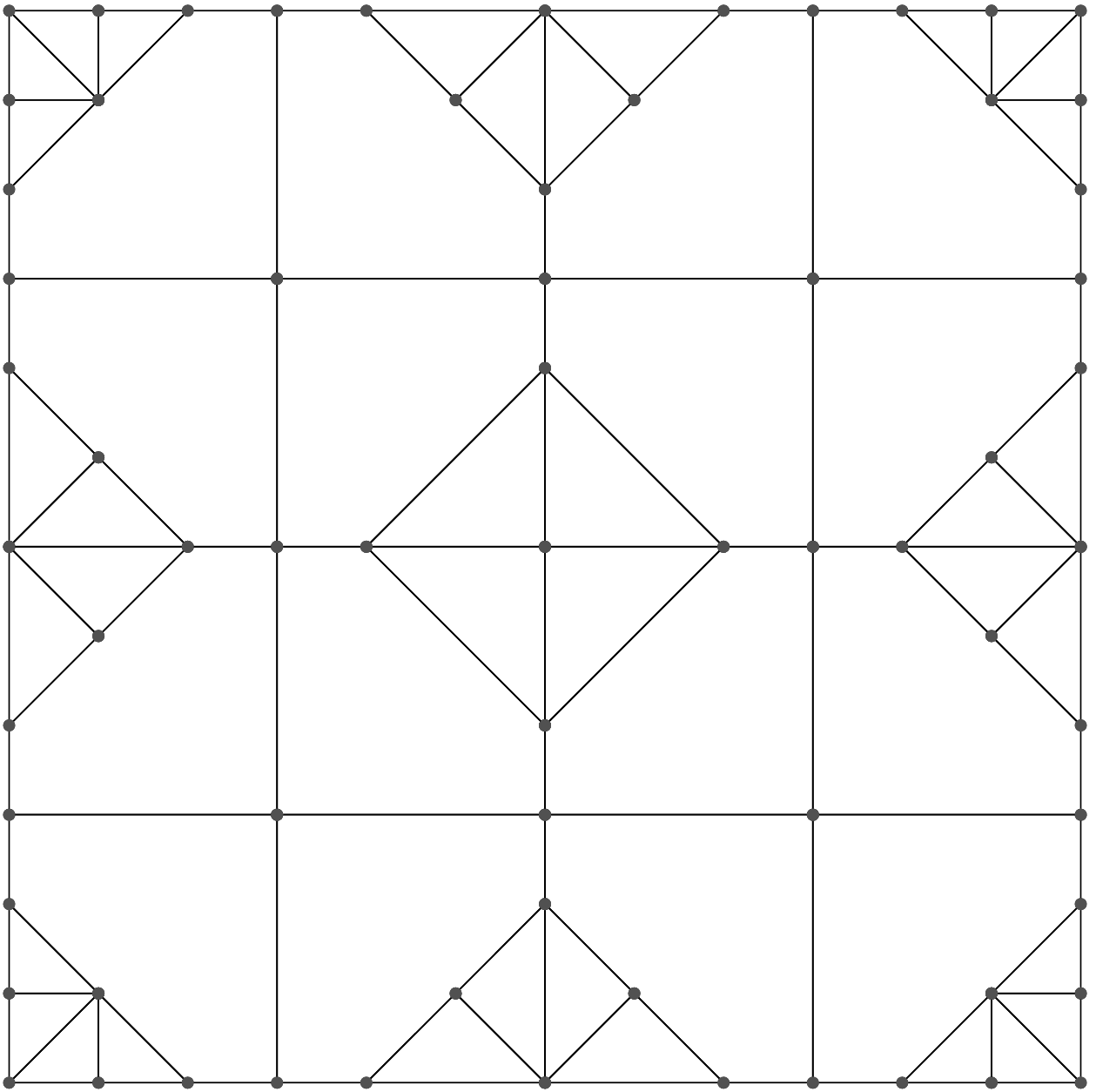}\hfill
 \includegraphics[trim=4.4cm 7.4cm 4.4cm 7.8cm, width=0.22\textwidth, clip]{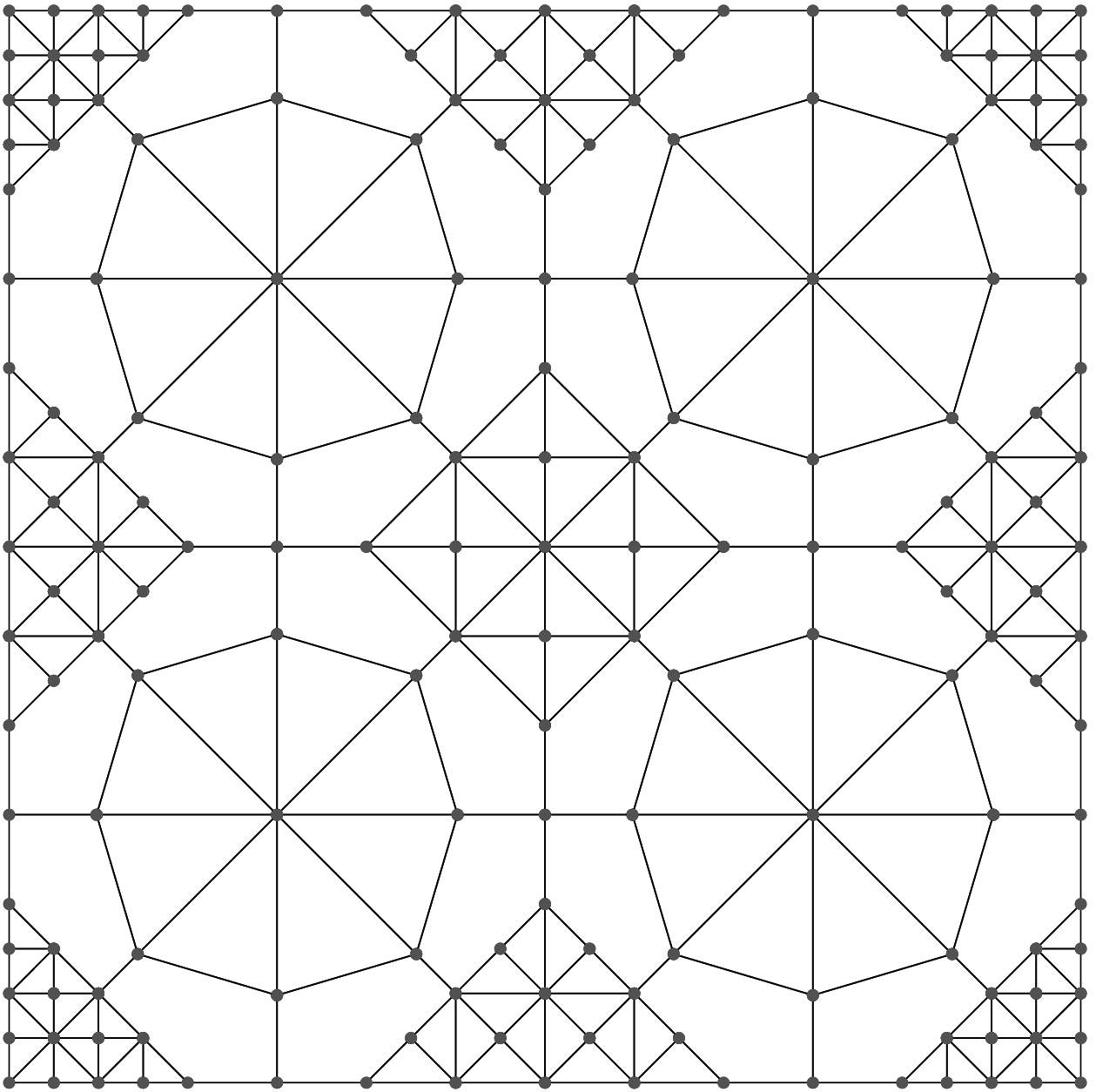}\hfill
 \includegraphics[trim=4.4cm 7.4cm 4.4cm 7.8cm, width=0.22\textwidth, clip]{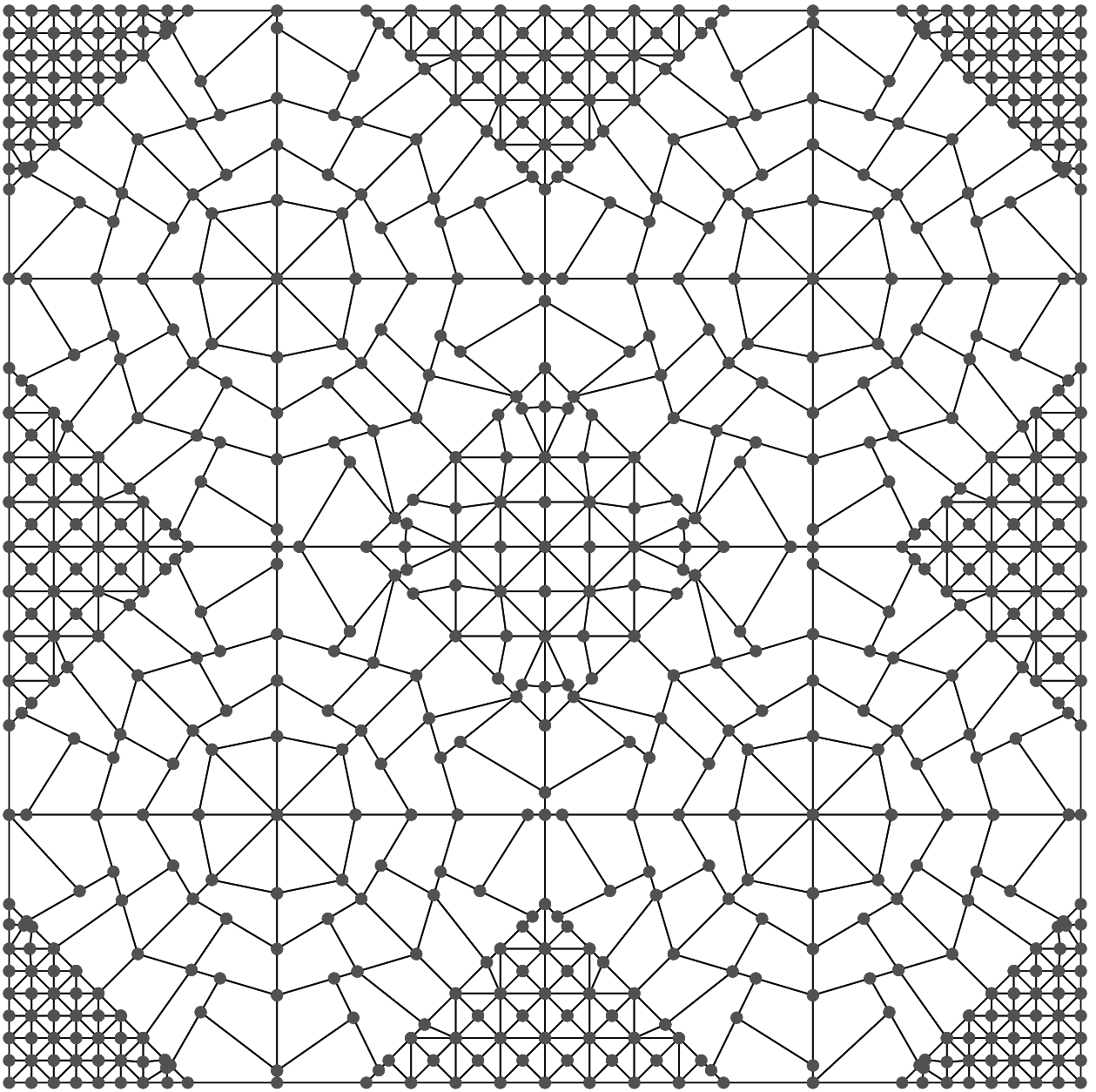}
 \caption{Initial mesh for Example~\ref{Example1} (left), and uniform refinements after $2,4$ and $6$ steps.}
 \label{fig:INT_Meshes}
\end{figure}
\begin{figure}[tbp]
 \centering
  \includegraphics[scale=0.8]{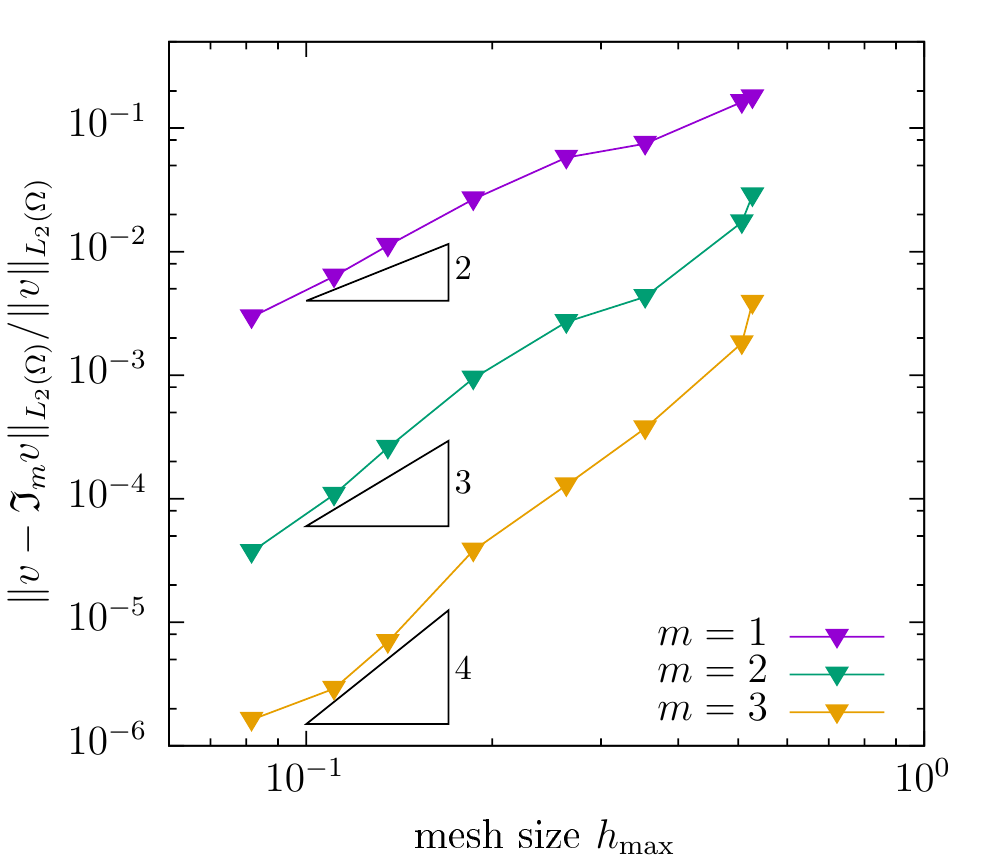}
  \includegraphics[scale=0.8]{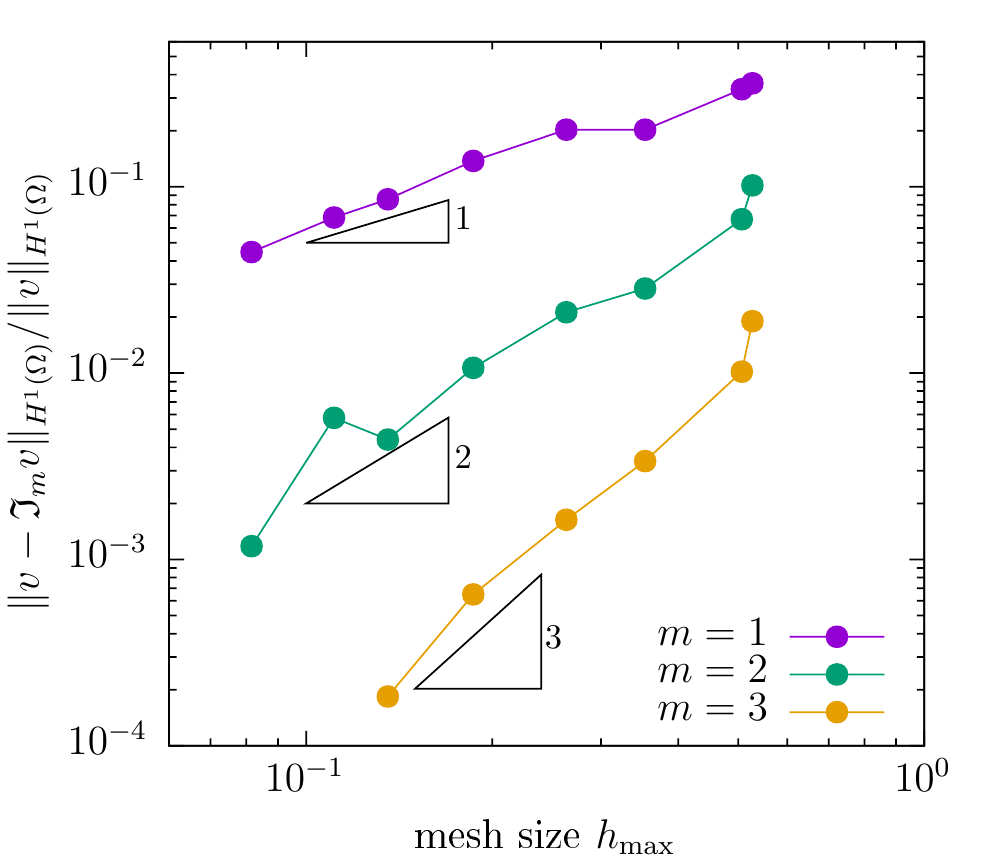}
 \caption{Relative interpolation error in $L_2$-  and
   $H^1$-norm for Example~\ref{Example1}, with $m=1,2,3$.}
 \label{fig:INT_Convergence}
\end{figure}
\end{example}

\begin{example}\label{Example2}
  In a second example, we interpolate the function
  $v(x) = r^{2/3}\sin(2(\theta-\pi/2)/3)$ on the L-shaped domain
  $\Omega=(-1,1)^2\setminus[0,1]^2$. This function exhibits the
  typical singularity at the reentrant corner.  For $m=1$, we compare
  the $L^2$-interpolation error for two families of meshes (see
  Figure~\ref{fig:L_Meshes}).  The $n^{th}$ mesh, $\cT_n$, of the
  first family consists of one L-shaped element,
  $(-1/3,1/3)^2\setminus[0,1/3]^2$, and $24n^2$ squares of size
  $(3n)^{-1}\times(3n)^{-1}$; $\cT_n$ has $(2n+1)(12n+1)+1$ vertices.
  The $n^{th}$ mesh, $\widehat\cT_n$, of the second family consists of
  congruent squares such that its number of vertices, which is of the
  form $(k+1)(3k+1)$, is as close as possible to $(2n+1)(12n+1)+1$.
  The dimensions of $V_1(\cT_n)$ and $V_1(\hat\cT_n)$ are clearly the
  number of vertices in the corresponding meshes.  For all square
  elements in either mesh, the local spaces are the bilinear
  functions.  If $K$ is the $L$-shaped element in $\cT_n$, then
  $\dim V_1(K)=6n+2$, and $V_1(K)$ contains functions having the
  correct singular behavior at the reentrant corner.

  We study the convergence of the relative interpolation error in
  $L_2(\Omega)$ with respect to the number of degrees of freedom (DoF)
  for both families of meshes. The optimal convergence behavior for the
  second family of meshes for an arbitrary smooth function is
  $\mathcal{O}(\mathrm{DoF}^{-1})$, but we neither expect or obtain
  that behavior for the given $v$, because of its singularity at the
  origin. In Table~\ref{table:conv:L_Meshes}, the relative interpolation error in
  the $L_2$-norm is given for the two sequences of meshes with
  comparable numbers of degrees of freedom. Furthermore, the numerical
  order of convergence (noc) is given. This is an estimate of the
  exponent $q$ in the error model $C\mathrm{DoF}^{-q}$.
  Since the function $v$ has a
  singularity, the convergence slows down for the standard bilinear
  elements on the uniform sequence of quadrilateral meshes. But, the
  optimal order of convergence is recovered for the uniform sequence
  with a fixed L-shaped element, because the local space associated
  with that element contains naturally functions with the correct singular
  behavior near the corner. 
\begin{figure}[tbp]
 \centering
 \includegraphics[width=1.2in]{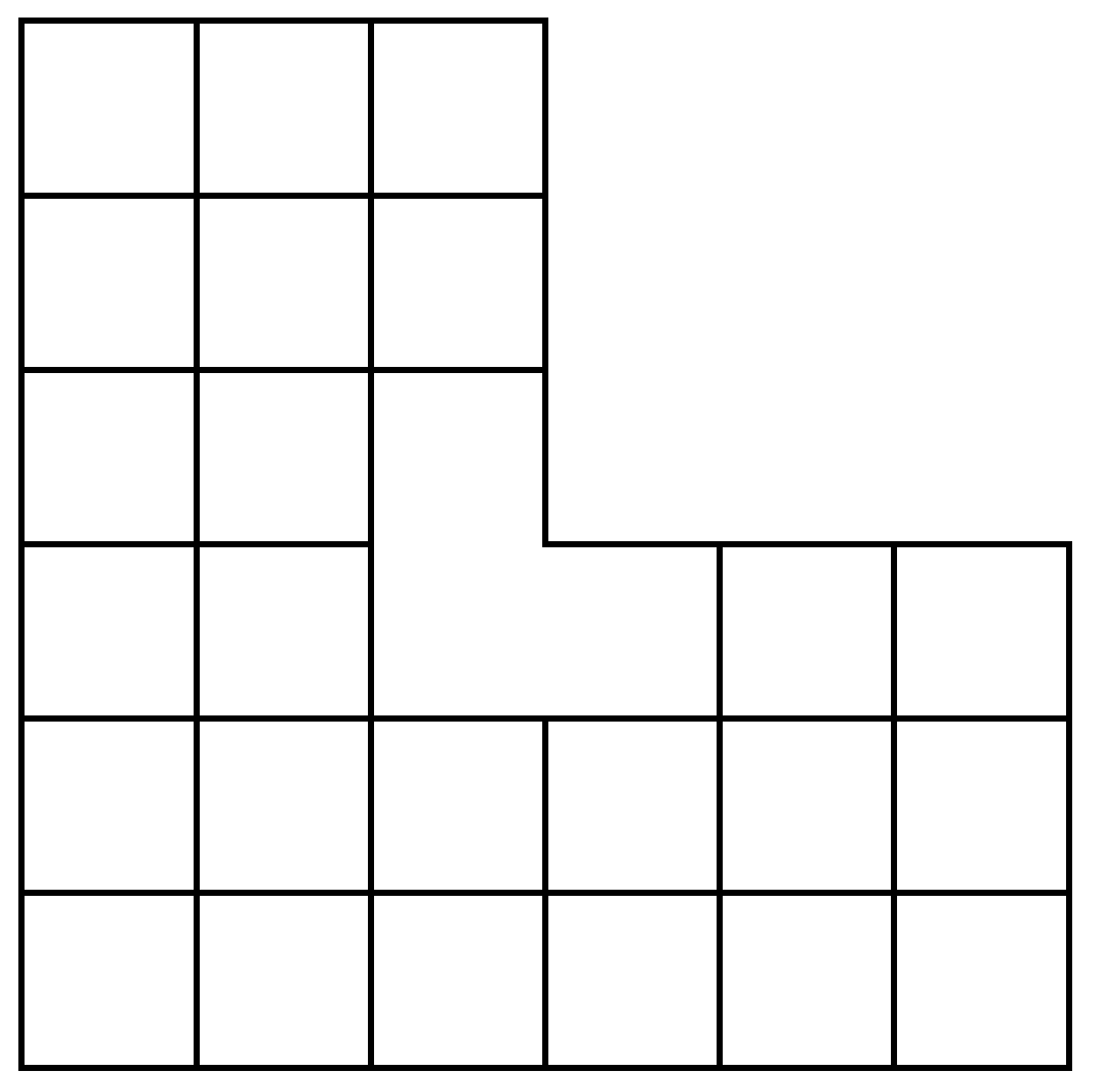}\hfil
 \includegraphics[width=1.2in]{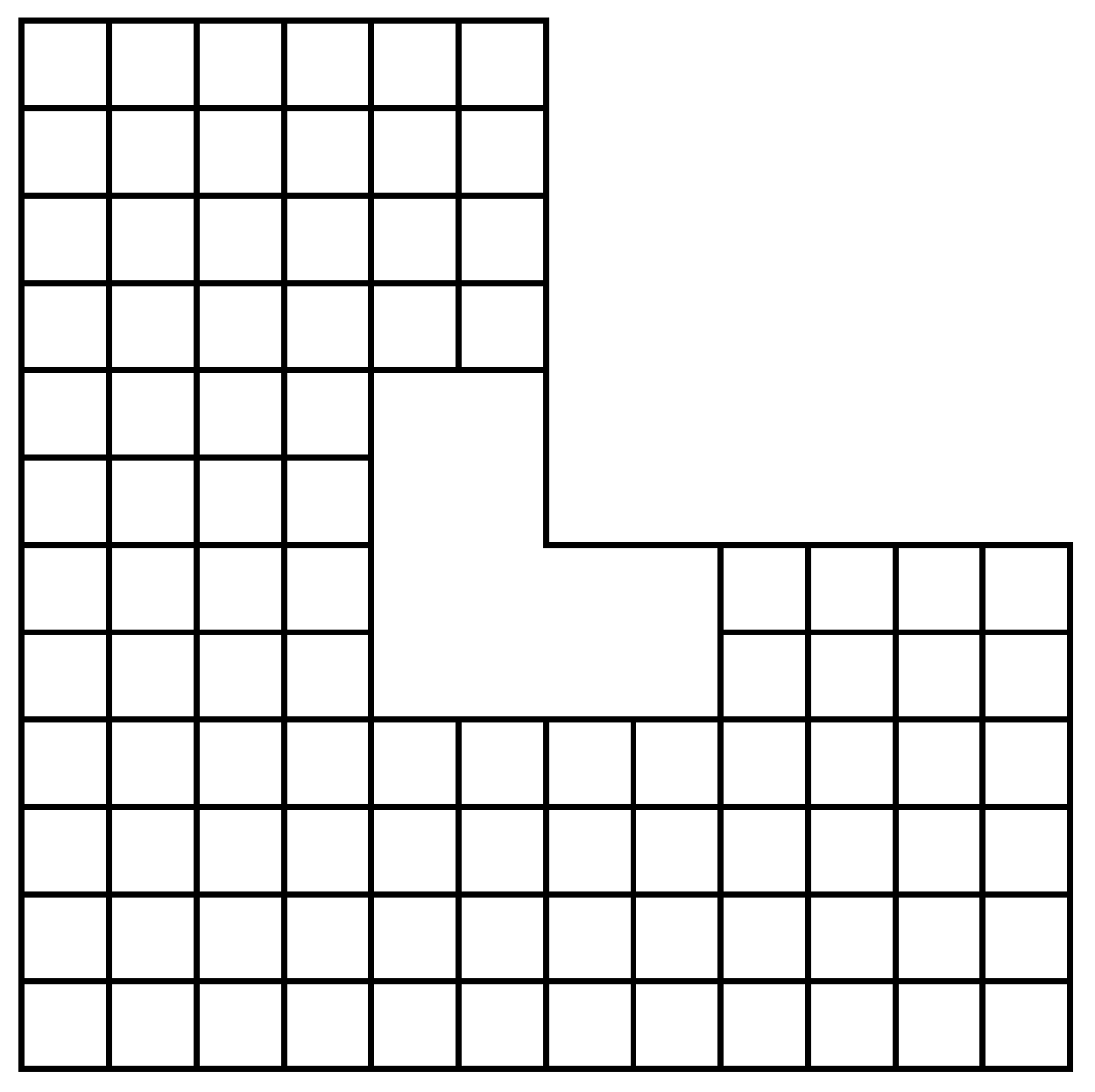}\hfil
 \includegraphics[width=1.2in]{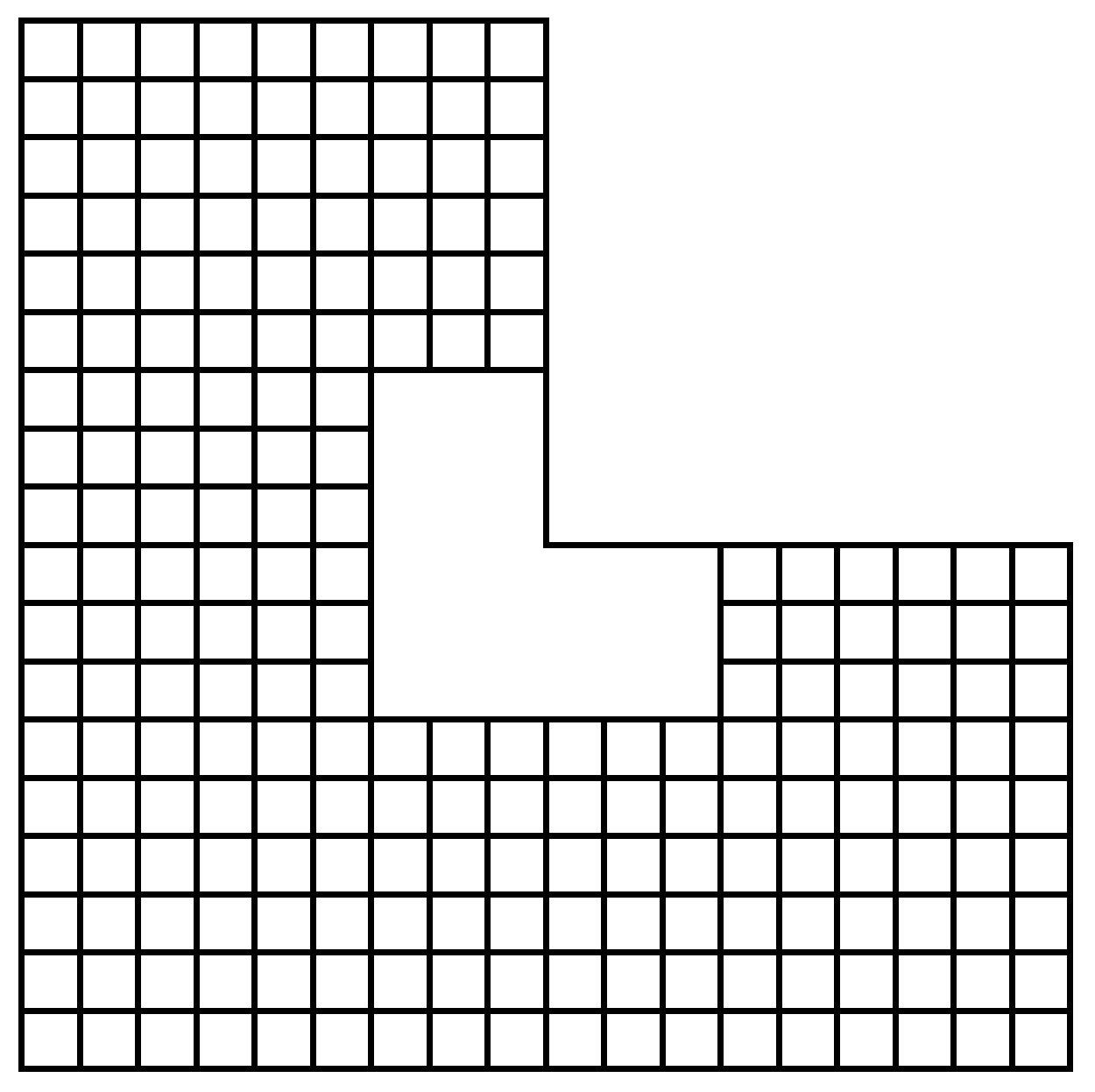}\hfil
 \includegraphics[width=1.2in]{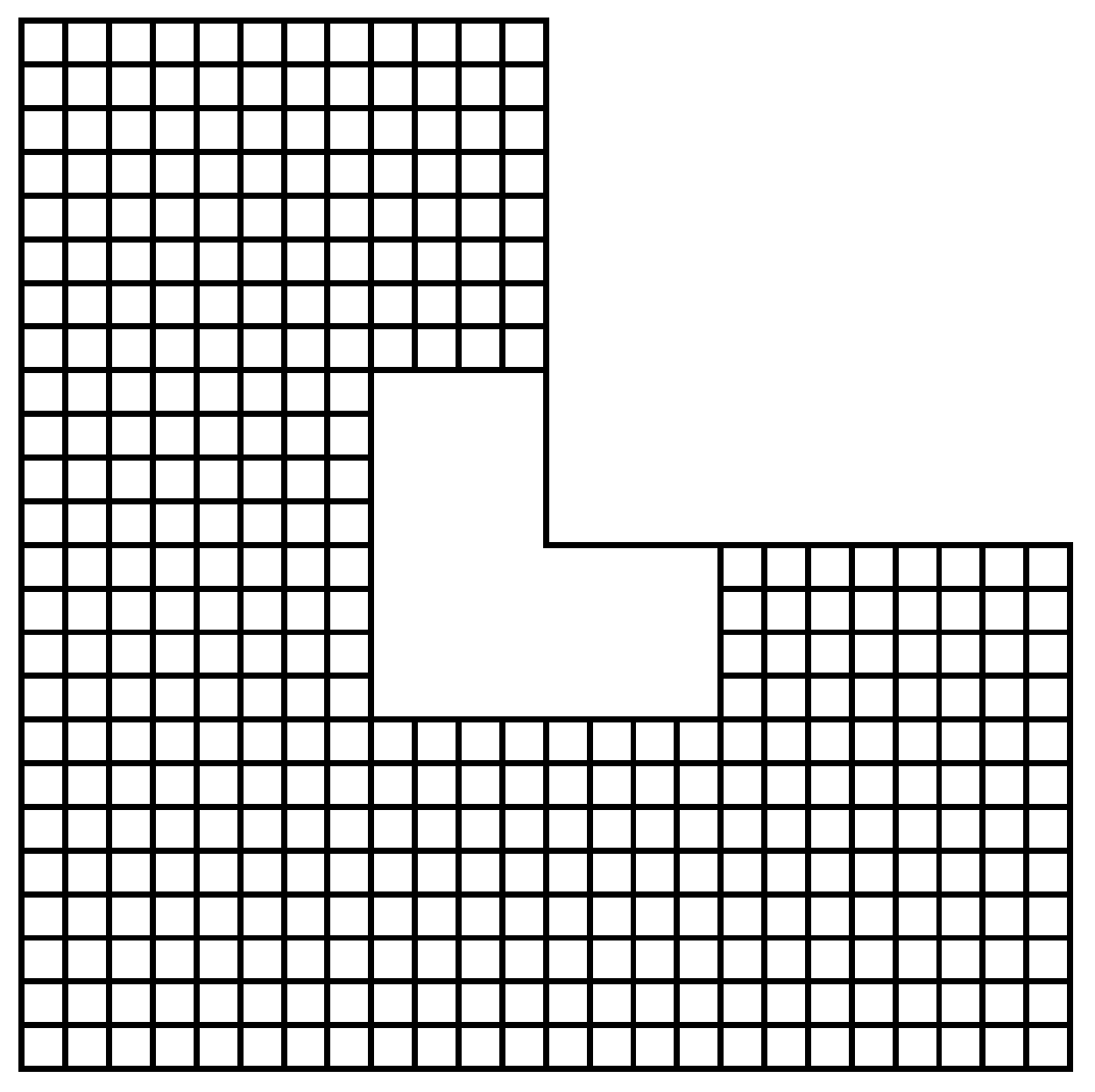}\\
\includegraphics[width=1.2in]{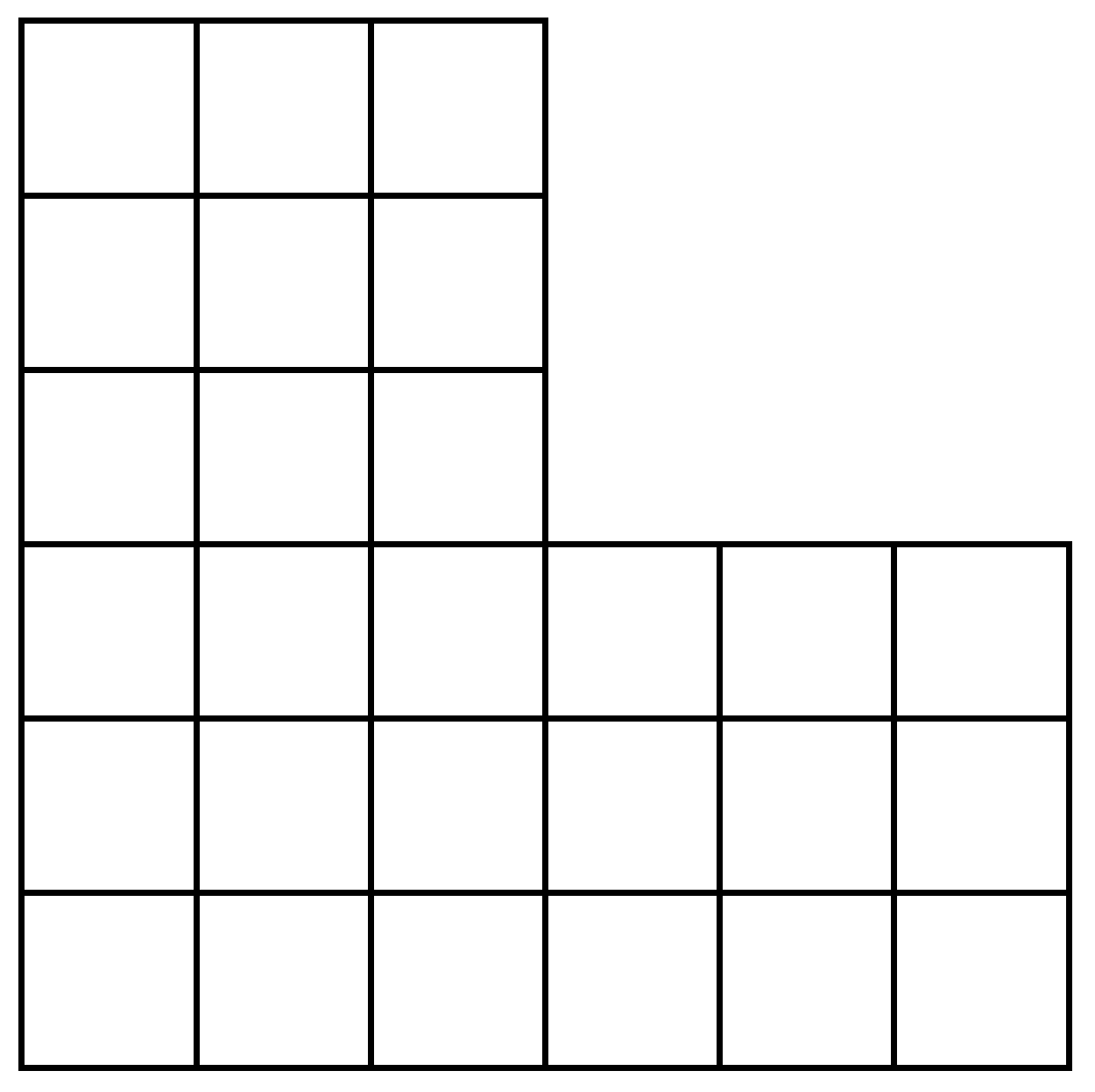}\hfil
 \includegraphics[width=1.2in]{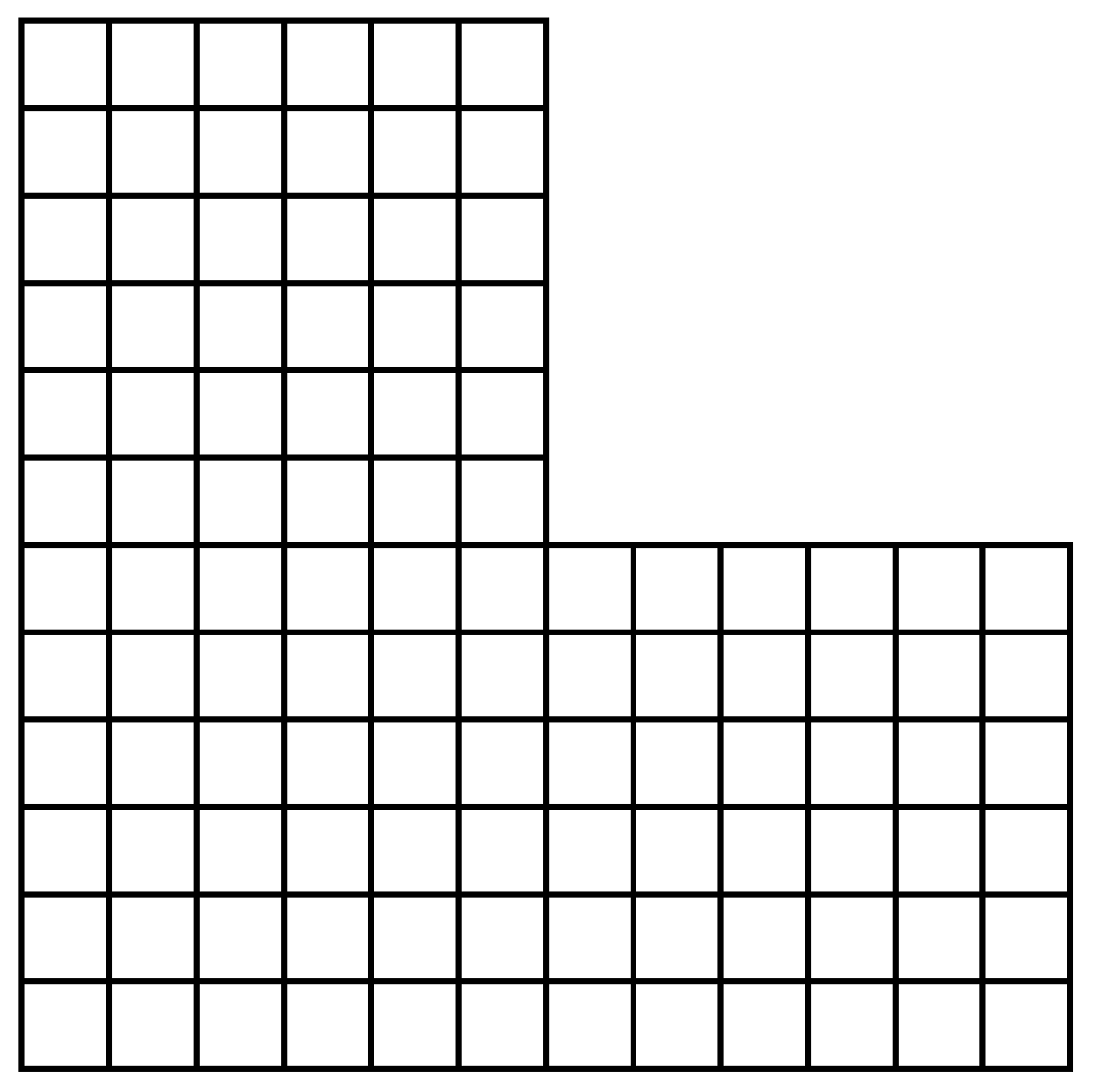}\hfil
 \includegraphics[width=1.2in]{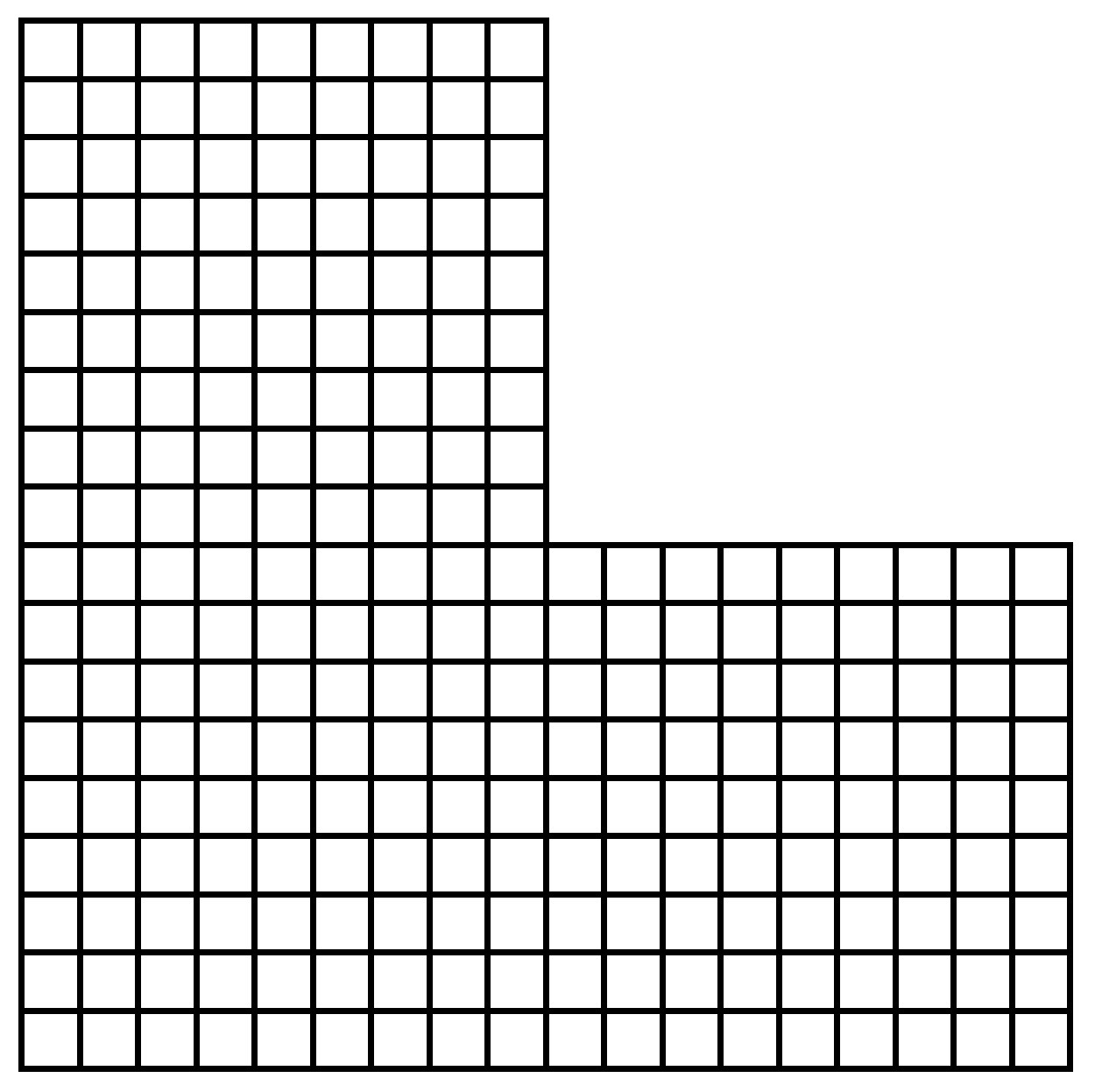}\hfil
 \includegraphics[width=1.2in]{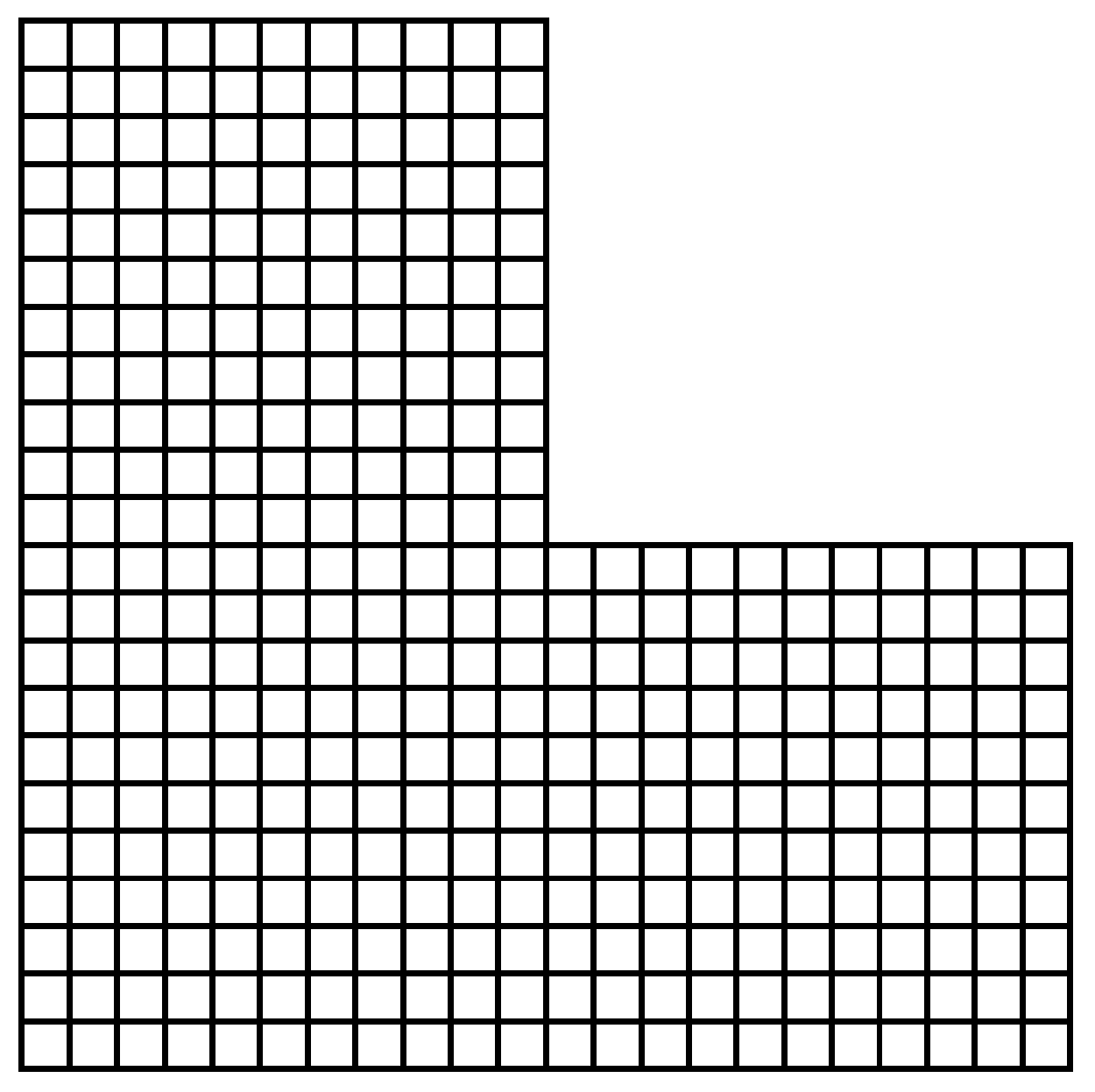}
 \caption{First four meshes of first family (top), and second family,
   for Example~\ref{Example2}. The fourth mesh in the first
   family has edge-length $1/12$ for each of its squares, and the fourth mesh in the second
   family has edge-length $1/11$ for each of its squares.}
 \label{fig:L_Meshes}
\end{figure}
\begin{table}
  \caption{Relative $L_2$-error (err) and numerical order of
    convergence (noc) with respect to the number of degrees of freedom
    (DoF) for first order interpolation ($m=1$) of  $v=r^{2/3}\sin(2(\theta-\pi/2)/3)$ in Example~\ref{Example2}.}
  \label{table:conv:L_Meshes}
\begin{tabular}{rcc|rcc}
  \hline
  \multicolumn{3}{c|}{First Family}& \multicolumn{3}{c}{Second Family}\\
  DoF & err & noc &  DoF & err & noc \\\hline
  $  40$ & 3.238e-03 &  --    & $  40$ & 1.263e-02 &  --    \\
  $ 126$ & 8.015e-04 & 1.22 & $ 133$ & 4.003e-03 & 0.96 \\
  $ 260$ & 3.549e-04 & 1.12 & $ 280$ & 2.042e-03 & 0.90 \\
  $ 442$ &1.989e-04 & 1.09 & $ 408$ & 1.463e-03 & 0.89 \\
  $ 672$ & 1.269e-04 & 1.07 & $ 833$ & 7.848e-04 & 0.87 \\
  $ 950$ & 8.789e-05 & 1.06 & $1045$ & 6.452e-04 & 0.86 \\
  $1276$ & 6.439e-05 & 1.05 & $1281$ & 5.415e-04 & 0.86 \\
  \hline
\end{tabular}
\end{table}
\end{example}

\section{Dirichlet Boundary Conditions, Curvilinear Elements}\label{Curvilinear}

As suggested in Remark~\ref{MixedSpace} and illustrated in
Example~\ref{NystromPacman}, the Nystr\"om approach for evaluating
functions that solve local Poisson problems readily accommodates
curvilinear elements and non-polynomial data.  As such, curved
boundaries or interior interfaces may be addressed more directly,
without resorting to polygonal approximations of these curves or
mappings from polygonal reference elements (e.g. isoparametric
elements, cf.~\cite{Scott1973,Lenoir1986,Bernardi1989}).  Although the
treatment of curved boundary and interior edges in our framework will
be investigated more thoroughly in later work, we here provide some
indication of how our approach may be used to address Dirichlet
boundary conditions on straight or curved edges, after first making a
few general remarks about necessary changes to the description of
$V_m(K)$ that must be made to accommodate curved edges.

We first remark that, if $K$ is a not a (straight-edged) polygon, the
definition of $V_m(K)$ must be either adjusted or properly
interpreted.  More specifically, if $e$ is a curved edge of $K$, the
defintion of $\PP_m(e)$, i.e. the polynomials of degree at most $m$ on
$e$, needs clarification.  One fairly natural approach is to define
$\PP_m(e)$ as the space of polynomials of degree at most $m$ \textit{
  with respect to arc length } on $e$.  We will call this approach the
\textit{Type 1} version of $\PP_m(e)$.  In this case,
$\dim \PP_m(e)=m+1$.  A potential drawback of this approach is that it does not
generally lead to the inclusion $\PP_m(K)\subset V_m(K)$, so we are
not guaranteed the approximation quality of $\PP_m(K)$.  A second
approach to defining $\PP_m(e)$ for a curved edge $e$ is to take it to
be the trace on $e$ of $\PP_m(\RR^2)$.  We will call this approach the
\textit{Type 2} version of $\PP_m(e)$. For the Type 2 version, we typically have $\dim \PP_m(e)=(m+2)(m+1)/2$, which leads
to a larger space $V_m(K)$, but yields the desired inclusion
$\PP_m(K)\subset V_m(K)$.  When $e$ is a straight edge, the Type 1 and
Type 2 versions of $\PP_m(e)$ yield the same space.
Note that one must revise the definition of
the degrees of freedom associated with a curved edge in the Type~2 case.
Basic differences between the two
approaches to defining $\PP_m(e)$ are illustrated in
Example~\ref{CircleSectorEx}.   A thorough investigation of these two
approaches, and of the practical and theoretical treatment of curved edges more generally, is a topic of subsequent work.

\begin{example}\label{CircleSectorEx}
  Let
  $K=K_h=\{x=(r\cos\theta,r\sin\theta):\;0<r<h\,,\,0<\theta<\alpha\pi\}$
  for some fixed $0<\alpha<2$, with two straight edges and one curved
  edge $e=\{x=(h\cos\theta,h\sin\theta):\;0\leq\theta\leq\alpha\pi\}$.
  We compare the two approaches to defining $\PP_m(e)$ in the case
  $m=1$.  For the Type 1 space, we have
  $P_1(e)=\mathrm{span}\{1,\theta\}$.  For the Type 2 space, we have
  $\PP_1(e)=\mathrm{span}\{1,\cos\theta,\sin\theta\}$.  Natural bases
  for these two approaches are plotted in
  Figure~\ref{CircleSectorFig}, with respect to $\theta$, for the
  choices $\alpha=1/2$ and $\alpha=3/2$.  The Type 1 basis functions
  are what one would expect, and undergo no qualitative changes as
  $\alpha$ varies.  The Type 2 basis functions were chosen with
  respect to the endpoints and midpoint of the edge, such that each
  basis function has the value one at one of these three points, and
  the value zero at the two others.  The qualitative behavior of Type
  2 basis functions clearly depends on $\alpha$ here.
\begin{figure}
\begin{center}
\includegraphics[width=2.5in]{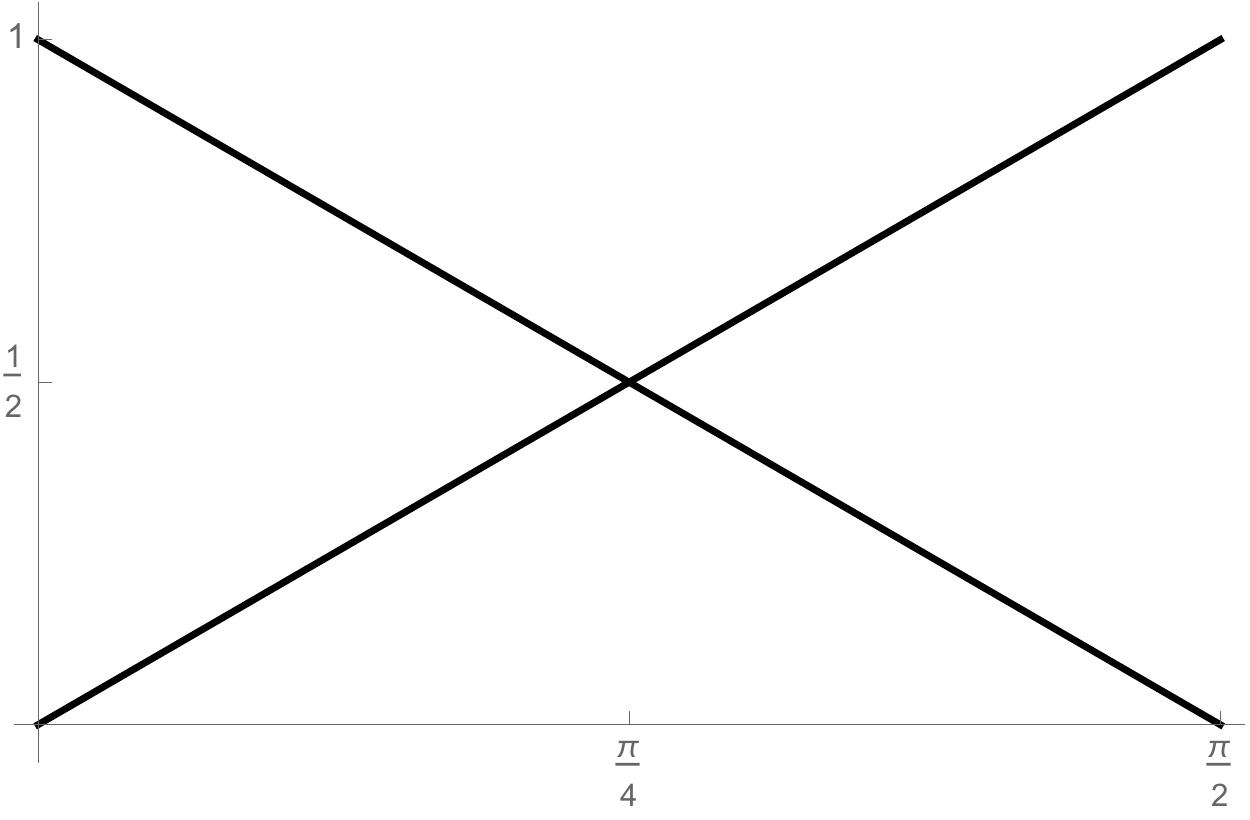}
\includegraphics[width=2.5in]{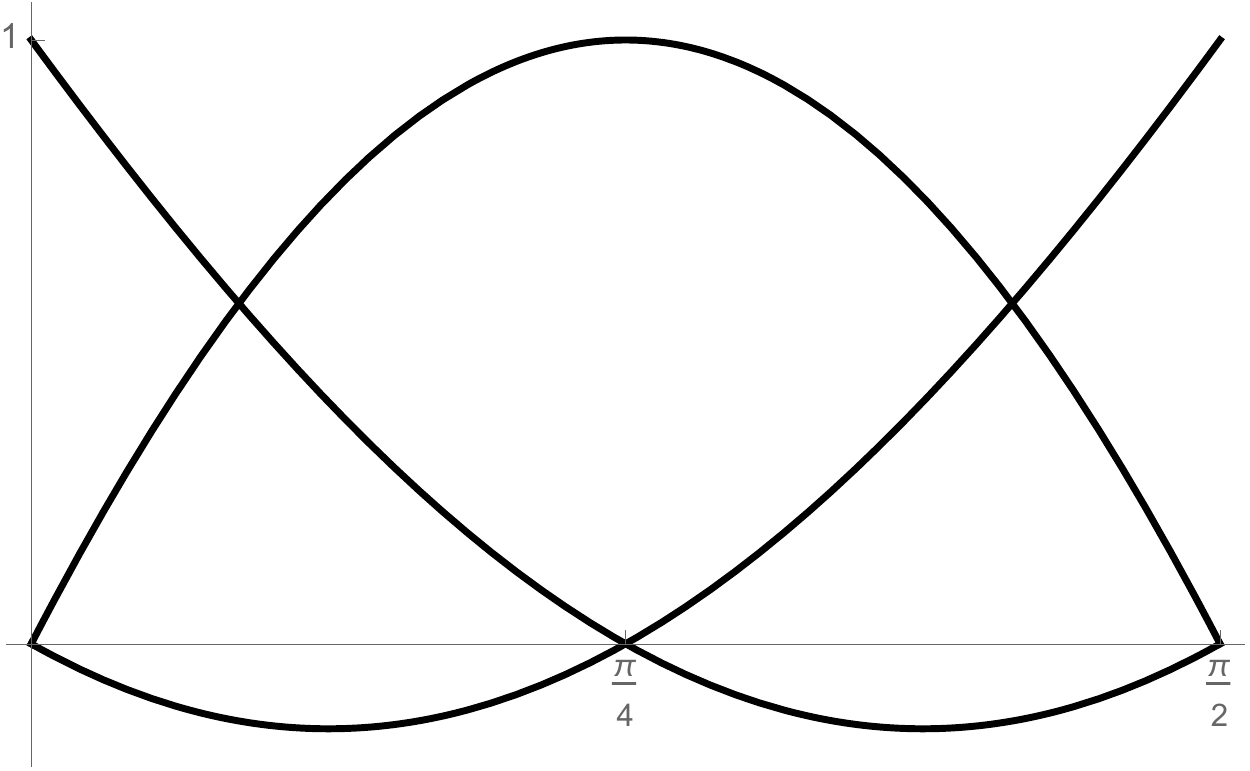}\\
\includegraphics[width=2.5in]{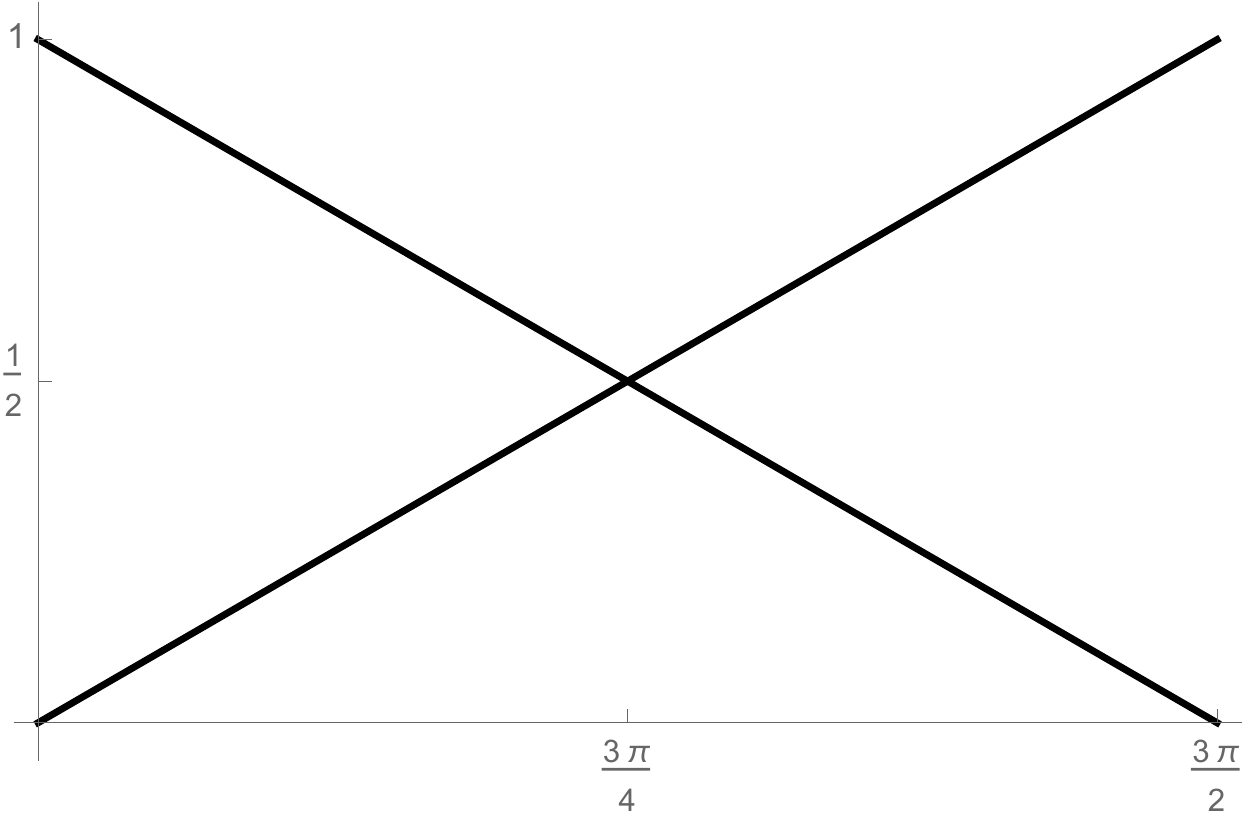}
\includegraphics[width=2.5in]{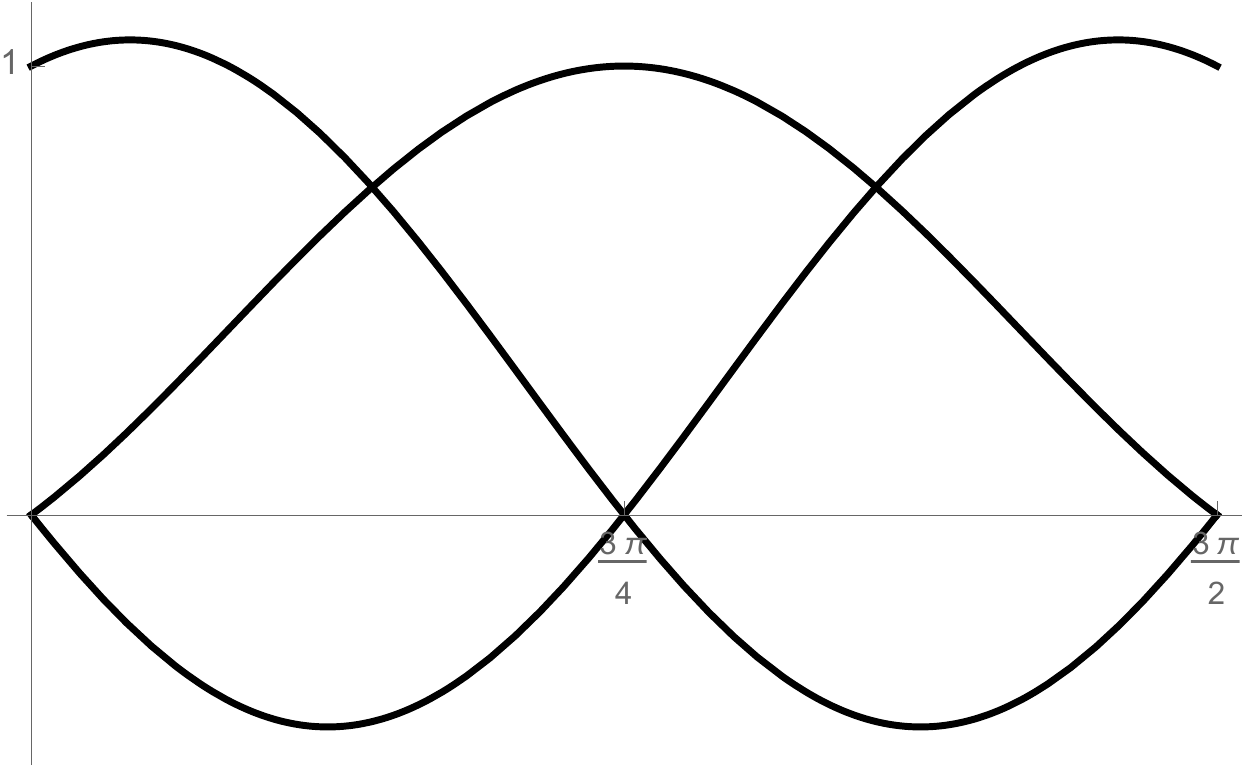}
\end{center}
\caption{\label{CircleSectorFig} Bases of the two types of spaces
  $\PP_1(e)$ for the circular arc in Example~\ref{CircleSectorEx}.
  First column corresponds to Type 1, and second column to Type 2.
  First row corresponds to $\alpha=1/2$, and second row to $\alpha=3/2$.}
\end{figure}
\end{example}

We now turn to the treatment of Dirichlet boundary conditions.
Let $K$ be curvilinear polygon such that $\partial K_D\doteq\partial
K\cap\partial\Omega\neq\emptyset$, and suppose we wish to 
prescribe boundary values $g$ that are continuous and
piecewise smooth on $\partial K_D$. At this stage, we assume that any curved edges are
contained in $\partial K_D$; so all interior edges are straight.
We take $g_D\in C(\partial K)$ to be $g$ on
$\partial K_D$, $0$ on edges not adjacent to $\partial K_D$, and
linear on edges adjacent to $\partial K_D$.
We employ the following local linear (and affine) spaces:
\begin{align}
V_m^K(K)&=\{v\in H^1_0(K):\;\Delta v\in\PP_{m-2}(K)\mbox{ in }K\}~,\\
V_m^{\partial K,0}(K)&=\{v\in H^1(K):\;\Delta v=0\mbox{ in }
K\;,\,v\in \PP_m^{0,D}(\partial K)\mbox{ on }\partial K \}~,\\
V_m^{\partial K,g}(K)&=v_g+V_m^{\partial K,0}(K)~,
\end{align}
where $\PP_m^{0,D}(\partial K)$ are those
functions in $C(\partial K)$ whose restriction to any edge
$e\not\subset\partial K_D$
is $\PP_m(e)$, and whose restriction to to any edge $e\subset\partial
K_D$ is $0$.   The function $v_g$ satisfies
\begin{align*}
\Delta v_g=0\mbox{ in }K\quad,\quad v=g_D\mbox{ on }\partial K~.
\end{align*}
The local affine space used in the global approximation is $V_{m,g,D}(K)=V_m^K(K)+V_m^{\partial K,g}(K)$.


\begin{example}\label{CurvedBoundaryExample}
We consider a single element 
\begin{align*}
K=K_h=\{(x_1,x_2)\in\RR^2:\;0\leq x_1\leq h\,,\,\sin(2\pi x_1/h)/4\leq
  x_2\leq h\}~,
\end{align*}
(see Figure~\ref{fig:CurvedMeshes1}).   We consider the convergence of
the interpolation error $v-\cI_1 v$ in $L^2(K)$ with respect to $h$
for $\cI_1 v\in V_{1,g,D}(K)$ and two different smooth functions $v$.
Since $m=1$, $\cI_1 v$ that agrees with $v$ on the curved edge,  
is equal to the linear interpolant of $v$ on each of the three
straight edges.  So computing $\cI_1 v$ requires the solution of a
single integral equation.  

If we had the inclusion $\PP_1(K)\subset
V_{1,g,D}(K)$, standard interpolation error estimates
(cf.~\cite[Theorem 1.103]{Ern2004}) would yield
\begin{align*}
\|v-\cI_1 v\|_{L^2(K)}\leq c h^2 |v|_{H^2(K)}
\end{align*}
for $v\in H^2(K)$.  Though we are not guaranteed, and will
typically not have, the inclusion $\PP_1(K)\subset V_{1,g,D}(K)$, we
still desire such quadratic convergence with respect to $h$ in practice.  The
experiments presented in Table~\ref{CurvedBoundaryInterpExper}
demonstrate this quadratic convergence in $h$ for two smooth functions
by considering ratios of successive errors as $h$ is halved.
\begin{figure}[tbp]
 \centering
 \includegraphics[width=1.5in]{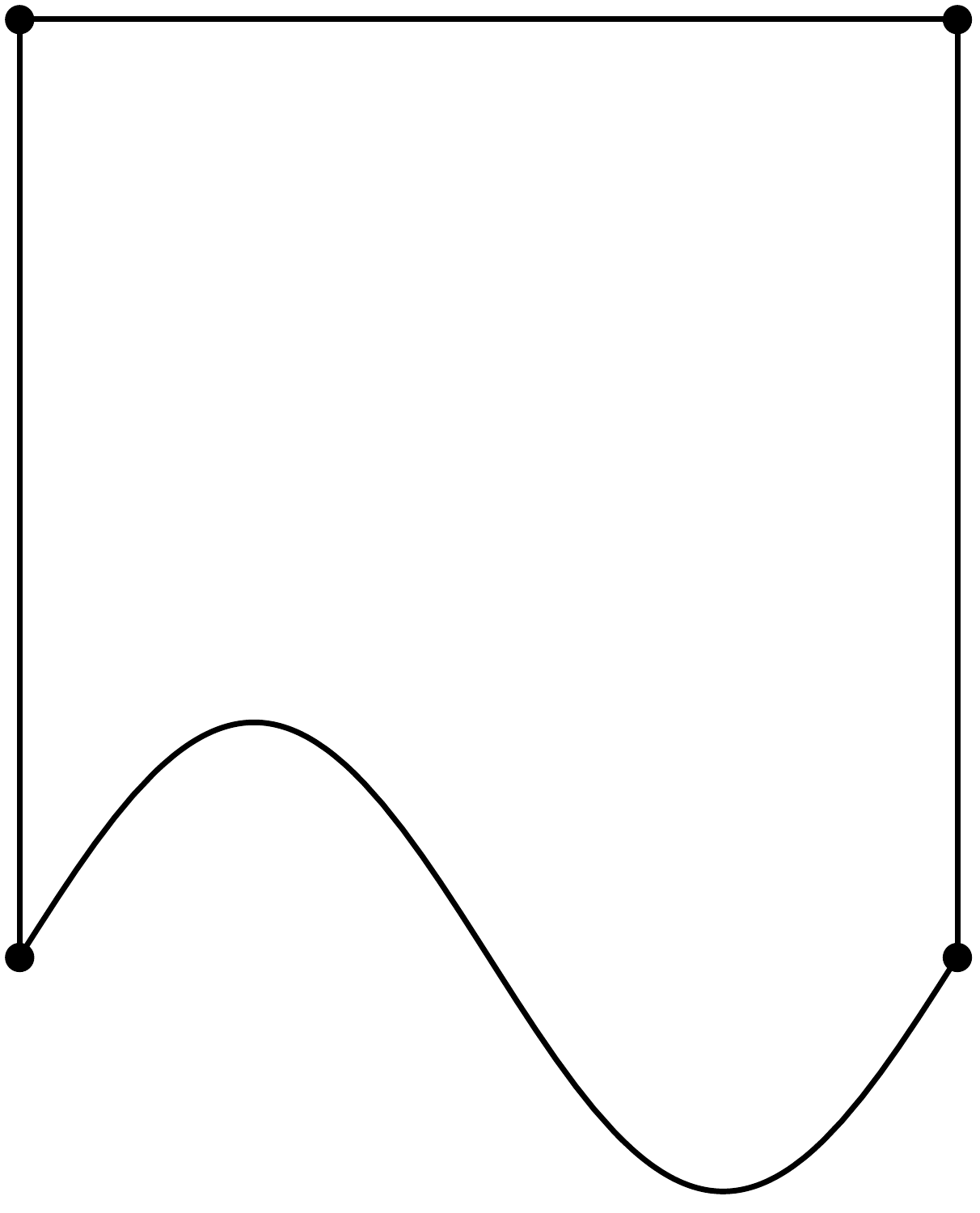}\hfil
  \caption{Element $K_h$ from Example~\ref{CurvedBoundaryExample}.}
 \label{fig:CurvedMeshes1}
\end{figure}
\begin{table}
\caption{\label{CurvedBoundaryInterpExper} $L^2$ interpolation error,
$\mathrm{error}=\|v-\cI_1 v\|_{L^2(K)}/ |v|_{H^2(K)}$, for $K=K_h$ from
Example~\ref{CurvedBoundaryExample}.}
\begin{tabular}{|c|cc|cc|}\hline
&\multicolumn{2}{c|}{$v=x^3-3xy^2+ 5 (x^2 -
  y^2)$}&\multicolumn{2}{c|}{$v=e^x+e^y$}\\
$h$&error&ratio &error&ratio\\\hline
$2^{-3}$&5.4199e-04&\rule[2.5pt]{0.8cm}{0.4pt}&1.5157e-03&\rule[2.5pt]{0.8cm}{0.4pt}\\
$2^{-4}$&1.3546e-04&4.0011&3.7899e-04&3.9994\\
$2^{-5}$&3.3856e-05&4.0010&9.4737e-05&4.0004\\
$2^{-6}$&8.4628e-06&4.0006&2.3682e-05&4.0004\\
$2^{-7}$&2.1155e-06&4.0003&5.9201e-06&4.0002\\
$2^{-8}$&5.2886e-07&4.0002&1.4800e-06&4.0001\\\hline
\end{tabular}
\end{table}
\end{example}

\section{Acknowledgements}
The work of JO on this paper was supported by the National Science
Foundation under Grant No. DMS-1414365.  AA thanks the NPDE-TCA for
supporting JO’s visit to IIT Kanpur that laid the foundation for this
paper.


\def\cprime{$'$}


\end{document}